\documentclass[11pt]{amsart}
\usepackage{graphicx}
\usepackage{epsfig}
\usepackage{epstopdf}
\usepackage{hyperref}
\usepackage{color}

\newtheorem{theorem}{Theorem}[section] 
\newtheorem{lemma}[theorem]{Lemma} 
\newtheorem{proposition}[theorem]{Proposition} 
\newtheorem{corollary}[theorem]{Corollary} 

\newcommand{\ignore}[1]{}

\begin{document}

\bibliographystyle{plain} 

\title{A lower bound technique for triangulations of simplotopes}

\author{Tyler Seacrest}
\address{The University of Montana Western\\
Dillon, MT 59725 \\
U.S.A.}
\email{tyler.seacrest@umwestern.edu}

\author{Francis Edward Su}
\thanks{The authors gratefully acknowledge the support of 
NSF Grants DMS-0301129, DMS-0701308, and DMS-1002938 (Su).}
\address{Department of Mathematics\\ Harvey Mudd College\\ Claremont, CA
91711\\ U.S.A.}
\email{su@math.hmc.edu}

                                                                             
\begin{abstract}
  Products of simplices, called \emph{simplotopes}, and their triangulations arise naturally in algorithmic applications in game theory and optimization.  We develop techniques to derive lower bounds for the size of simplicial covers and triangulations of simplotopes, including those with interior vertices. We establish that a minimal triangulation of a product of two simplices is given by a 
  \emph{vertex triangulation}, i.e., one without interior vertices.
  For products of more than two simplices, we produce bounds for products of segments and triangles.  Aside from cubes, these are the first known lower bounds for
  triangulations of simplotopes with three or more factors, and our techniques suggest extensions to products of other kinds of simplices.  We also construct a minimal triangulation of size 10 for the product of a triangle and a square using our lower bound.
\end{abstract} 

\maketitle

\section{Introduction}
A classical problem in discrete geometry is to determine the size of a
minimal triangulation of a given polytope.
For instance, a polytope that has received considerable attention is
the $d$-dimensional cube; see e.g., 
\cite{blsu05, cott82, huan96, sall82, smit00} for upper and
lower bounds on the size of many kinds of minimal decompositions of the cube.
Minimal triangulations also serve a practical purpose as well; they
can be used in simplicial algorithms for finding fixed points (e.g.,
see \cite{todd76, yang96}) as well as for economic applications \cite{su99} 
and applications to GIS~\cite{ohst16},
since smaller triangulations lead to  more efficient algorithms.
By a \emph{triangulation}, we mean a decomposition of a polytope $P$ into simplices that meet face-to-face, and we allow the vertex set of the triangulation to include more points than just the vertices of $P$.  
(Triangulations that only use vertices of $P$ will be called \emph{vertex triangulations}.)

In this paper, we study minimal triangulations of 
{\em simplotopes}, which are products of simplices \cite{freu86}.  
The $d$-cube is thus a special kind of simplotope, the product of $d$ 
$1$-dimensional simplices (segments).  
Simplotopes are of special
interest in economics, since in a non-cooperative $n$-person game, the space of strategies is the product of simplices (one for each player)---
and finding a Nash-equilibrium is equivalent to finding a fixed point of a function on this space, e.g., \cite{lata82}.
Simplotopes and their triangulations also appear in algebraic geometry \cite{arbi07, babi98, stur96} and optimization {\cite{dkos09}.
Orden and Santos \cite{orsa03}
used an ``efficient" $38$-simplex triangulation of a simplotope--- the
product of a $3$-cube and a triangle--- to construct triangulations of
arbitrarily high-dimensional cubes with few simplices; however, it
should be noted the concept of an efficient triangulation is
different from (though related to) the minimal triangulation.
The idea of triangulating products of polytopes by using triangulations of simplotopes as building blocks is present also in \cite{haim91, sant00}.

Aside from cubes, very little is known about minimal triangulations of 
simplotopes, especially for triangulations that are not vertex triangulations.   
It is well-known (e.g., ~\cite{orsa03}) that the product
of two simplices of dimensions $a$ and $b$ must be triangulated with exactly ${a + b \choose a}$
simplices if it is a vertex triangulation;  one of our results in this paper is that a vertex triangulation is indeed minimal over all triangulations.
Work has been done to enumerate the many different vertex triangulations of such a product~\cite{taii97}. 
The space of triangulations of the product of two simplices has a rich structure.  For example, Dyck paths lead to an important class of such triangulations~\cite{ceps15}.   This case is also important in combinatorial commutative algebra (see \cite{cort14} for example) and is a ``building block" for more general polytopal constructions (see Section 6.2 of \cite{ders10}).
DeLoera, Rambau and Santos \cite{ders10} give a recent survey of the enumerative and 
structural properties of triangulations of polytopes.

\section{Summary of Results}
Let $\Delta^d$ denote the standard $d$-dimensional {\em simplex}, 
the convex hull of the $d+1$ standard basis vectors in $\mathbb{R}^{d+1}$.  
Thus ${\bf x}=(x_1,...,x_{d+1})$ in $\Delta^d$ satisfies $x_i\geq 0, \sum_i x_i =1$.
A \emph{simplotope} is the product of simplices: 
$\Delta^{c_1} \times \cdots \times \Delta^{c_n}$.
%Thus, each point in $\Delta^{c_1} \times \cdots \times \Delta^{c_n}$ corresponds to a selection of one point in each factor $\Delta^{c_i}$.
%For example, a point in $\Pi(2, 2)$ corresponds to a choice of a pair of points, one from each triangle.
%See Figure \ref{tricrosstri}. 
% At times it will be helpful to have a shorthand for simplotopes with many repeated factors. 
We use the shorthand $\Pi(a_1, \ldots, a_n)$ to represent the product
$$
(\Delta^1)^{a_1} \times (\Delta^2)^{a_2} \times \cdots \times (\Delta^n)^{a_n}
$$
which is the product of $a_1$ segments, $a_2$ triangles, $a_3$ tetrahedrons, etc.  
Much of the paper will focus on $\Pi(s, t)$, which is the product of $s$ segments and $t$ triangles.

Two kinds of simplotopes will receive special attention: 
(i) products of two simplices of any dimension and 
(ii) arbitrary products of segments and triangles.

We shall obtain lower bounds for the size of a minimal triangulation by studying the associated concept of a \emph{cover}.  
Given a $d$-dimensional polytope $P$, a collection of $d$-simplices is a 
\emph{(vertex) cover} of $P$ if the union of the simplices is $P$ and the
vertices of each simplex are vertices of $P$.  (All covers in this paper will be vertex covers, so we shall just refer to them as covers.) Thus a vertex triangulation is a special kind of cover of $P$.
The following result provides the key to connect the study of covers to the study of triangulations, and it 
may be surprising, in light of the fact that a general triangulation of $P$ is not necessarily a cover of $P$.

\begin{theorem}[Bliss-Su]
\label{bliss-su}
Let $P$ be a polytope.  
Let $C(P)$ be the size of the minimal cover of $P$, using only simplices spanned by the vertices of $P$.
Let $T(P)$ be the size of a minimal triangulation of $P$, possibly using vertices that are not vertices of $P$.
Then
$$C(P) \leq T(P).$$
\end{theorem}

\begin{table}[h]                                                                

\bigskip

 \tiny                                                                       
 \begin{center}                                                                 
 \begin{tabular}{|c|llllllll|}\hline                                             
 \setlength{\unitlength}{1mm}                                                   
 \begin{picture}(2,6)                                                           
 \thinlines                                                                     
 \put(2, 3){t}                                                                  
 \put(-1.5, 0){s}                                                               
 \put(-2.9, 5.8){\line(1, -1){6.8}}                                             
 \end{picture}                                                                  
& 0 & 1 & 2 & 3 & 4 & 5 & 6 & 7 \\ \hline                                           
$0$ & $1 $ & $1$ & $6$ & $50$ & $423$ & $4240$ & $50179$ & $732543$ \\
$1$ & $1$ & $3$ & $20$ & $163$ & $1523$ & $16467$ & $232398$ & $3267672$ \\
$2$ & $2$ & $9$ & $68$ & $612$ & $6048$ & $74586$ & $989216$ & \\
$3$ & $5$ & $33$ & $256$ & $2499$ & $28333$ & $355024$ & $5118845$ & \\
$4$ & $16$ & $136$ & $1054$ & $10784$ & $120878$ & $1629902$ & &\\
$5$ & $60$ & $532$ & $4552$ & $48713$ & $576338$ & $9090179$ & &\\
$6$ & $250$ & $2203$ & $19416$ & $216057$ & $3096452$ & & &\\
$7$ & $1117$ & $8897$ & $92047$ & $1134649$ & $16362555$ & & &\\
$8$ & $4680$ & $44740$ & $450047$ & $6001487$ & & & &\\
$9$ & $21384$ & $218063$ & $2392586$ & $36139265$ & & & &\\
$10$ & $95708$ & $1144311$ & $13708102$ & & & & &\\
$11$ & $516465$ & $5853664$ & $75322333$ & & & & &\\
$12$ & $2906455$ & $33135045$ & & & & & &\\
$13$ & $16372399$ & $197669956$ & & & & & &\\
$14$ & $91944719$ & & & & & & & \\ 
$15$ & $522902357$ & & & & & & & \\  \hline                  
\end{tabular}                                                                   
\end{center}                                                                    
\bigskip
%\normalsize 
\caption{
\label{results-table}
Our linear program produces these lower bounds for the number of
simplices needed to cover (and hence to triangulate) $\Pi(s,t)$, the product of $s$ segments and $t$ triangles,
for various
values of $s$ and $t$.  
Compare these with other known results from Table \ref{comparison-table}.
}
\end{table}

\begin{table}[h]                                                                

\bigskip

 \small                                      
 \begin{center}                                                                 
 \begin{tabular}{|c|lll|}\hline                                             
 \setlength{\unitlength}{1mm}                                                   
 \begin{picture}(2,6)                                                           
 \thinlines                                                                     
 \put(2, 3){t}                                                                  
 \put(-1.5, 0){s}                                                               
 \put(-2.9, 5.8){\line(1, -1){6.8}}                                             
 \end{picture}                                                                  
& 0 (cubes) & 1 & 2\\ \hline                                           
$0$ & $C,T=1$ & $C, T=1$ & $C, T=6$ \\
%$1$ & $C,T=1$ & $C, T=3$ & $C \geq 20, T \leq 26^{(1)}$  \\
$1$ & $C,T=1$ & $C, T=3$ & $T \leq 26^{(*)}$  \\
$2$ & $C,T=2$ & $C, T =10^{(\dagger)}$ & \\
% $3$ & $C,T=5$ & $C \geq 33, T \leq 38$\textsuperscript{\cite{orsa03}} & \\
$3$ & $C,T=5$ & $T \leq 38$\textsuperscript{\cite{orsa03}} & \\
$4$ & $C,T=16$\textsuperscript{\cite{blsu05}} & & \\
$5$ & $C \geq 60$\textsuperscript{\cite{blsu05}}, $T = 67$\textsuperscript{\cite{huan96}}  & &  \\
$6$ & $C \geq 252$\textsuperscript{\cite{blsu05}}, $T = 308$\textsuperscript{\cite{huan96}}  & &  \\
$7$ & $C \geq 1143$\textsuperscript{\cite{blsu05}}, $T = 1493$\textsuperscript{\cite{huan96}}  & &  \\
$8$ & $C \geq 5104$\textsuperscript{\cite{blsu05}}, $T \leq 11944$\textsuperscript{\cite{orsa03}}  & &  \\
$9$ & $C \geq 22616$\textsuperscript{\cite{blsu05}}, $T \leq 173015$\textsuperscript{\cite{sall82}}  & & \\
$10$ & $C \geq 98183$\textsuperscript{\cite{blsu05}}, $T \leq 1728604$\textsuperscript{\cite{sall82}} & &  \\        \hline  
\end{tabular}                                                                   
\end{center}                                                                    
\bigskip
%\normalsize 
\caption{
\label{comparison-table}
Other known results for the covering number $C$ and the triangulation number $T$ of various simplotopes.   Note that $C \leq T$ always holds by the Bliss-Su Theorem.  For cubes, lower bounds from~\cite{blsu05} are slightly better than our lower bounds from Table \ref{results-table} due to techniques specific to cubes.
These are taken from the citations in superscript, computations using TOPCOM$^{(*)}$ 
(see~\cite{ramb02} to learn about TOPCOM), 
and results in Section~\ref{tri-cross-sq} of this paper$^{(\dagger)}$.  }

\end{table}

Bliss and Su proved this result in \cite{blsu05} by considering a piecewise linear map taking vertices 
of a triangulation to the vertices of a cover, 
using a Sperner labelling of the vertices, and showing that this map has degree 1.  
They used it to obtain bounds for minimal triangulations of cubes; 
we develop new techniques to extend their ideas to find bounds for minimal covers and triangulations of simplotopes.  

We can now explain how this paper is organized.
In Section \ref{background}, we establish terminology and background on simplotopes---including coordinate representations, volume considerations, and the standard triangulation---that serve as a foundation for the rest of the paper.  

Then in Section \ref{product-two-simplices} we study products of two simplices and, via Theorem \ref{bliss-su}, we observe that any vertex triangulation (indeed the standard triangulation) is a minimal triangulation (Theorem \ref{product-of-2-simps}).  Previously it was known that any two vertex triangulations must have the same number of simplices, but had not been known if adding extra vertices could reduce the size of a triangulation.

We then develop techniques in Sections \ref{tools} and \ref{recurrence-relation} that extend methods for cubes developed by Bliss and Su \cite{blsu05} to study arbitrary products of segments and triangles.  
Although this isn't the most general simplotope, it will be apparent that with further work our techniques could be used to study simplotopes that are products of other kinds of simplices, and when possible we state our results for general simplotopes.  Our analysis yields linear programs that arise from considerations of covering exterior faces and exploiting the product structure of these polytopes.  

%These tools allow us in Section~\ref{recurrence-relation} to develop a recurrence relation 
%to be used in a linear program to give the lower bounds found 
%in the table in Figure~\ref{results-table}.  

Our results are summarized in Table \ref{results-table}, in which we
provide lower bounds for triangulations of $\Pi(s,t)$ for several small dimensions.  
Aside from cubes, these are the first known lower bounds for triangulations of simplotopes with three or more factors.  Computationally, we obtain the bounds by computing the values of a certain recursive function and then solving a linear program based on those values.

As an example, the efficient $38$-simplex triangulation of $\Pi(3,1)$ by Orden and Santos \cite{orsa03} (for vertex triangulations) compares favorably to our lower bound of 33 in Table \ref{results-table} (for triangulations that allow extra vertices).

We list our bounds in Table \ref{results-table} up to dimension $15$.  At this dimension, the linear program is still easy to compute, but the increasing gap between upper and lower bounds makes the results beyond this dimension less interesting. Moreover, as the dimension grows beyond this, the computational bottleneck arises first with the recursive function, not the linear program.

We end the paper in Section \ref{tri-cross-sq} by studying in detail $\Pi(2, 1)$, the product of a square and triangle.  We construct a size $10$ triangulation and prove it is minimal by using our lower bound of $9$ from Table \ref{results-table} together with additional geometric arguments.

\section{Background}
\label{background}

\subsection*{Coordinate Representation}
%There are two different coordinate representations that we will find convenient for representing simplotopes.

The {\em standard coordinate representation} 
expresses each point $\mathbf{v}$ of $\Delta^{c_1} \times \cdots \times \Delta^{c_n}$ as
% a $(c_1+...+c_n + n)$-tuple in $\mathbb{R}^{c_1+...+c_n + n}$:
$$
\mathbf{v} = (\mathbf{x}^1; \mathbf{x}^2; \ldots; \mathbf{x}^n)
$$
where each $\mathbf{x}^i$ is a point in $\Delta^{c_i}$.
We write $\mathbf{x}^i = (x_1^i, \ldots, x_{c_i+1}^i)$ and say the
coordinates within each $\mathbf{x}^i$ are in the same \emph{factor}.
We separate the coordinates of different factors by semi-colons.
Note that for each point $\mathbf{v}$, 
we have that
$x_j^i \geq 0$ for all $i$ and $j$ and $x_1^i + \cdots + x_{c_i+1}^i = 1$ for
all $i$.  The latter relations imply that although we 
represent the simplotope as an object in $\mathbb{R}^{c_1+...+c_n + n}$, 
it has dimension $c_1+...+c_n$.  Additionally, if all $x_j^i$ are integers (0 or 1), 
then $\mathbf{v}$ is a vertex of $\Delta^{c_1} \times \cdots \times \Delta^{c_n}$.

For example, the vector $(0.2, 0.3, 0.5; 0.1, 0.3, 0.6)$ represents a point in 
the interior of $\Delta^2 \times \Delta^2$ in the standard coordinate representation.
Notice the two triplet factors in this vector sum to one.  
Another point in $\Delta^2 \times \Delta^2$ is the vertex $(0, 1, 0; 0, 0, 1)$.  See Figure \ref{tricrosstri}.

\begin{figure}[h]
\begin{center}
\includegraphics[scale=.6]{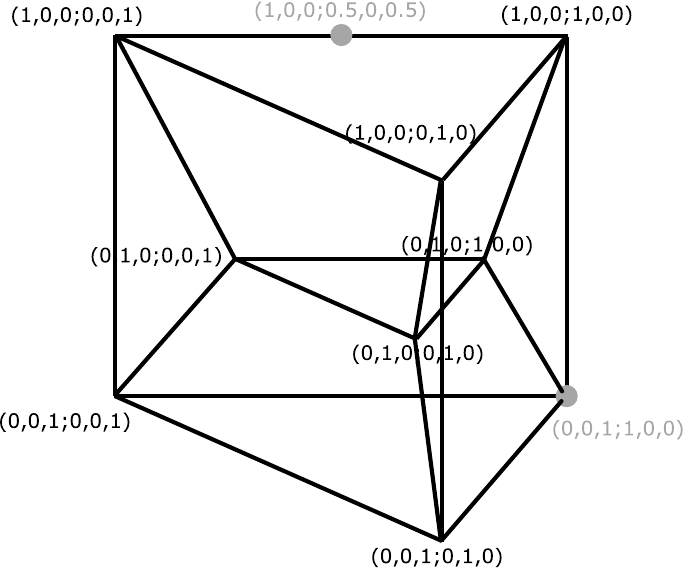}
\end{center}
\caption{
\label{tricrosstri}A Schlegel diagram of $\Delta^2 \times \Delta^2$, a 4-dimensional polytope, labeled by standard coordinates in $\mathbb{R}^6$.
The two grey points lie on the same $2$-face.}
\end{figure}

%For some purposes we will want to use a \emph{reduced coordinate
%  representation} that expresses points of $\Delta^{c_1} \times \cdots \times \Delta^{c_n}$
%by vectors in $\mathbb{R}^{c_1+...+c_n}$.  This is similar to
%the standard representation, except for each $i$ exactly one of
%$x_1^i$, \ldots, $x_{c_i+1}^i$ is removed, i.e., in each factor, one
%of its coordinates is removed.  Which coordinates are removed can be
%specified by picking a vertex of the simplotope and ``forgetting'' all
%the coordinates that are non-zero for that vertex.  We call this the
%{\em reduced coordinate representation of $\Delta^{c_1} \times \cdots \times \Delta^{c_n}$ with
%  respect to a vertex $v$} and call the process {\em reduction with
%  respect to $v$}.
%
%As an example, consider reducing the standard representation of
%$\Pi(2, 2)=\Pi^*_{0,2}$ with respect to the vertex $v=(0,0,1; 0,1,0)$.  Then
%the point $(0.2, 0.3, 0.5; 0.1, 0.3, 0.6)$ is expressed in 
%reduced coordinates (with respect to $v$) 
%as $(0.2, 0.3; 0.1, 0.6)$.  In this
%reduction, the last coordinate in the first factor and the second
%coordinate in the last factor are removed, because those
%correspond to the non-zero coordinates of $v$.
%
%Note that no information
%is lost because the removed coordinates can always be recovered by
%remembering that the standard coordinates in each factor must sum to $1$.

Any collection of points $\mathbf{v}_1, \ldots, \mathbf{v}_k$ in standard 
coordinates can be represented in matrix form as
$$
M(\mathbf{v}_1, \ldots, \mathbf{v}_k) = \left [ \begin{array}{c}
    \mathbf{v}_1 \\ \vdots \\ \mathbf{v}_k \end{array} \right ]
$$
in which the rows of the matrix are the given points in standard coordinates.
%Similarly, we let 
%$M_\mathbf{v}(\mathbf{v}_1, \ldots, \mathbf{v}_k)$ denote the
%\emph{reduced matrix} whose rows are the given points in reduced
%coordinates with respect to a vertex $\mathbf{v}$.  There is a linear transformation (a projection) taking $M$ to $M_v$.

%The simplotope $\Delta^{c_1} \times \cdots \times \Delta^{c_n}$ is defined by the half-spaces
%$x_j^i \geq 0$, one half-space for each $i, j$ combination,
%intersected with the hyperplanes $\sum_{j=1}^{c_i} x_j^i = 1$, one
%hyperplane for each $i$.
The simplotope $\Delta^{c_1} \times \cdots \times \Delta^{c_n}$ is defined by the intersection of the half-spaces $x_j^i \geq 0$ and the hyperplanes $\sum_{j=1}^{c_i} x_j^i = 1$.
Note that a \emph{$k$-face} is the intersection of the simplotope and
($c_1 + ... +c_n - k$) of the hyperplanes $x_j^i = 0$.  See Figure \ref{standard} for an example.
The following lemma easily follows from noting that each column of the matrix representation corresponds with a coordinate $x_j^i$.

\begin{figure}[h]
\begin{center}
\includegraphics[scale=.75]{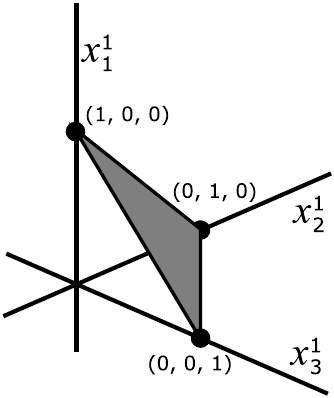}
\end{center}

\caption{\label{standard}An equilateral triangle in the
  standard coordinate representation.  Note that it is defined by
  the plane $x_1^1 + x_2^1 + x_3^1 = 1$, as well as the half-spaces
  $x_1^1 \geq 0$, $x_2^1 \geq 0$, and $x_3^1 \geq 0$.}
\end{figure}

\begin{lemma}
\label{exterior}
The set of points, $\mathbf{v}_1, \ldots, \mathbf{v}_m$ in 
$\Delta^{c_1} \times \cdots \times \Delta^{c_n}$ lie on the same $k$-face 
if and only if their standard
matrix representation $M(\mathbf{v}_1, \ldots, \mathbf{v}_m)$ has at
least $(c_1 + \cdots + c_n - k)$ columns consisting of only zeros. 
\end{lemma}

For example, with $c_1 = 2$ and $c_2 = 2$ (see Figure
\ref{tricrosstri}), the two points $(1, 0, 0; 0.5, 0, 0.5)$ and 
$(0, 0, 1; 1, 0, 0)$ have two coordinates, $x_2^1$ and $x_2^2$, 
that are zero in both points.  
Therefore they lie in the same $2 + 2 - 2 = 2$-face.

%\begin{proof}
%By definition, the points $\mathbf{v}_1, \ldots, \mathbf{v}_m$ lie on the
%same $k$-face of the simplotope $\Delta^{c_1} \times \cdots \times \Delta^{c_n}$ if and only if
%each point lies on the intersection of $(c_1 + \cdots + c_n- k)$ hyperplanes of the form $x_j^i = 0$.
%For each such hyperplane $x_j^i = 0$, the matrix representation $M(\mathbf{v}_1, \ldots, \mathbf{v}_m)$ for
%these points will have corresponding $ij$-th column equal to zero. 
%\end{proof}

\subsection*{Exterior Faces of Simplices}
For simplotopes, we say a $n$-simplex of a triangulation 
of polytope $P$ 
has an \emph{exterior} $j$-face if the simplex has $j+1$
vertices in the same $j$-face of $P$.

\subsection*{The Standard Triangulation}

There is a standard triangulation of any simplotope $\Delta^{c_1} \times \cdots \times \Delta^{c_n}$ of size
\begin{equation}
\label{multi-choose}
\frac{(c_1 + c_2 + \cdots + c_n)!}{c_1! c_2! \cdots c_n!}.
\end{equation}
This is also called the \emph{staircase triangulation}~\cite{ders10}.  To demonstrate, we introduce a new
coordinate system that will allow for a very simple permutation
description of the simplices in the standard triangulation.  
Note that this coordinate representation will be different than the
standard coordinate system defined earlier.
We can represent a point $\mathbf{w}$ in $\Delta^{c_1} \times \cdots \times \Delta^{c_n}$ by
a vector
$$
\mathbf{w} = (\mathbf{y}^1; \mathbf{y}^2; \ldots; \mathbf{y}^n)
$$
where $\mathbf{y}^i = (y_1^i, \ldots, y_{c_i}^i)$, each $y_j^i \in [0, 1]$,  
and 
\begin{equation}
\label{simplo-ineq}
y_{j-1}^i \geq y_j^i
\end{equation}
for all $i \in \{1, \ldots, n\}$ and $j \in \{2, \ldots, c_i\}$.
Hence the coordinates within each factor, from left to right, are non-increasing.  
For example, in $\Delta^2 \times \Delta^2$  we have that
$y_1^1 \geq y_2^1 \geq y_3^1$, and that $y_1^2 \geq y_2^2 \geq y_3^2$.
However, the relative sizes of the coordinates of $\mathbf{y}^1$ and
$\mathbf{y}^2$ are unrelated.
Note that the restrictions on each factor $\mathbf{y}^i$ define a
point in the simplex that is the convex hull of the $c_i+1$ points
\begin{eqnarray*}
\mathbf{u}_0 & = & (0, 0, 0, \ldots, 0) \\
\mathbf{u}_1 & = & (1, 0, 0, \ldots, 0) \\
& \vdots & \\
\mathbf{u}_{c_i} & = & (1, 1, 1, \ldots, 1).
\end{eqnarray*}
Hence, $\mathbf{w}$ defines a point of a simplotope.
Note that if all $y_j^i$ are integers, then $\mathbf{w}$ is a vertex of
$\Delta^{c_1} \times \cdots \times \Delta^{c_n}$.

By placing further restrictions on the coordinates 
%of the standard representation of the simplotope, 
we may obtain a subdivision that is
a triangulation.  Consider any ordering of {\em all} 
the coordinates $y_j^i$ that
is consistent with the inequalities in (\ref{simplo-ineq}).  
An example of such an ordering if $c_1 = 2$ and $c_2 = 2$ is
$$
y_1^2 \geq y_2^2 \geq y_1^1 \geq y_3^2 \geq y_2^1 \geq y_3^1.
$$
One may check that: (i) each such ordering defines a simplex, 
(ii) every point in the simplotope satisfies at least one such
ordering, and is therefore in (at least) one of the simplices, 
(iii) these simplices meet face-to-face.

Therefore these simplices form a triangulation of the simplotope,
called the {\em standard triangulation}.
The number of simplices in the triangulation is equal to the number of
ways of arranging the $y_j^i$'s, subject to the prior arrangement of
the coordinates within each factor being in non-increasing order.
This amounts to choosing positions in the ordering 
for the coordinates within each factor, and this is given by 
the multi-choose expression (\ref{multi-choose}).

The standard triangulation is the largest possible vertex triangulation of a simplotope; because the simplices of this triangulation have the smallest possible volume.

\subsection*{Class}
One way to obtain bounds for the number of simplices required to
triangulate $\Delta^{c_1} \times \cdots \times \Delta^{c_n}$ is to use volume estimates.
For example, the volume of the smallest $k$-simplex 
with vertices at lattice points is $\frac{1}{k!}$, and therefore 
the size of any triangulation with vertices at lattice points 
cannot be larger than 
%the volume of the polytope divided by $\frac{1}{k!}$, or
$k!$ times the volume of the polytope.

For convenience, we use a kind of normalized volume so that volumes of lattice point simplices 
are integers.  We define the {\em class} of a $d$-dimensional set in $\mathbb{R}^d$ to be the volume of that set multiplied by $d!$ (this concept is also referred to as {\em lattice volume} in the literature.)
Starting with a simplex $\alpha$ whose vertices are rows of the matrix $M$, choose any vertex $\mathbf{v}$ and create a new matrix $M_\mathbf{v}$ by removing every column where $\mathbf{v}$ contains a $1$, and also removing the row that corresponds to $\mathbf{v}$.  We define the \emph{class} of $\alpha$ to be $|\det(M_\mathbf{v})|$, the absolute value of the determinant of $M_\mathbf{v}$.  Note that this determinant is well-known to be $d!$ times the volume of the simplex~\cite{stein66}.

If $\alpha$ is a $k$-dimensional simplex embedded in $\mathbb{R}^d$ for $k < d$, the can define this notion of class similarly:  simply remove all the columns of zeros from $M_{\mathbf{v}}$ before taking the determinant.  

%Then if $\alpha$ is a simplex whose vertices are specified (in
%a reduced coordinate system) by the rows of a matrix $M_\mathbf{v}$,
%the {\em class} of $\alpha$ will be $|\det[1|M_\mathbf{v}]|$, which denotes the absolute value of
%determinant of the matrix formed by augmenting $M_\mathbf{v}$ 
%with a column of ones in front.  So in a vertex triangulation or a cover, the class of any simplex must be an integer, because it is the determinant of 
%an integer matrix $[1|M_\mathbf{v}]$ filled with $1$'s and $0$'s.
%The class can be obtained from the $\mathbb{R}^{c_1+...+c_n}$-volume of the the simplex  
%spanned by the rows of the reduced matrix $M_\mathbf{v}$ by multiplying by $(c_1 + \cdots + c_n)!$ ~\cite{blsu05}.
Although our definition of class appears to depend on the choice the vertex $\mathbf{v}$ of the simplex, since class corresponds to volume, this choice does not matter:
\begin{lemma}
Let $\alpha$ be a simplex of a vertex triangulation of a simplotope.
Then the class of $\alpha$ is independent of the choice of reduction $\mathbf{v}$ used to compute it.
\end{lemma}
%
%\begin{proof}
%This follows easily from noting that $[M_\mathbf{v_i}]$ and $[M_\mathbf{v_j}]$ are related by elementary column operations that do not involve scaling by numbers other than $1$ or $-1$, since entries in certain columns are replaced by the first column minus the sum of the columns of one factor.
%\end{proof}

\section{Products of Two Simplices}
\label{product-two-simplices}

These ideas can be used to demonstrate the well-known
result that for the product of two simplices, $\Delta^a \times \Delta^b$,
every vertex triangulation has the same number of simplices (e.g., see \cite{ders10, orsa03}).  From the standard matrix
representation and Lemma \ref{exterior}, one can verify that any simplex $\sigma$ using vertices of $\Delta^a \times \Delta^b$ must have an exterior facet---
the standard matrix $M$ of $\sigma$ has  $a+b+1$ rows each with two $1$'s in them, and $a+b+2$ columns, so 
some column must have no more than one $1$ in it.  
Removing the row of that $1$, if needed, produces an exterior facet with a single $0$ coordinate (there is at most one $0$ coordinate because the $\sigma$ is non-degenerate).
In Proposition~\ref{facetsame}, we will show that if
a simplex has an exterior facet, then the class of the simplex is
equal to the class of that facet.  But that facet is a simplex in a facet of the simplotope, hence in a product
of two simplices of lower total dimension.  
Then induction can be used to show that every 
simplex in $\Delta^a \times \Delta^b$ has class one.  Hence, to
cover the entire volume of $\Delta^a \times \Delta^b$, one must use a number of
simplices equal to the number of simplices in the standard
triangulation.  Hence every vertex triangulation is a minimal triangulation.

Since our argument above is a covering argument, it actually shows a stronger result via Theorem \ref{bliss-su} that has not been proved before:

\begin{theorem}
\label{product-of-2-simps}
A minimal triangulation of a product of two simplices cannot be smaller than a vertex triangulation (e.g., the standard triangulation).
\end{theorem}

%For instance, the introduction of extra vertices inside the $7$-dimensional polytope $\Delta^3 \times \Delta^4$ 
%will not allow for a smaller triangulation than the standard triangulation.

Thus the product of two simplices is uninteresting in the sense that the standard triangulation is a minimal triangulation.
This is not true for the product of three or more simplices, because
in it there exist simplices of class larger than one.
For the product of three segments (the $3$-cube) the minimal
triangulation is $5$ while the standard triangulation is $6$, because
the former uses a class 2 simplex, which has no exterior facet.  In the
product of three triangles, it is possible to have simplices of
class 4.
% (determined by a brute force computer search).

\section{Tools for Covering Simplotopes}
\label{tools}

Our goal is to determine a lower bound on the number of simplices
required in a cover of $\Pi(s,t)$, the product of $s$ segments and $t$ triangles.
By Theorem \ref{bliss-su}, such a lower bound is
also a lower bound for triangulations of $\Pi(s,t)$, even when
interior vertices are allowed.  In this section, we will develop several tools needed to do so, when possible stating results for general simplotopes.  Recall $\Pi(a_1, \ldots, a_n)$ to be the product of $a_1$ segments, $a_2$ triangles, $a_3$ tetrahedrons, and so on.

Notice that every face of $\Pi(s,t)$ is also a simplotope that is
the product of segments and triangles.  
Consider a $\Pi(s',t')$ face 
of $\Pi(s,t)$.
Note that while $t \geq t'$, it is not necessarily true that $s \geq
s'$.  This is because when one coordinate of a triangle factor is
fixed at $0$, then possible values for the remaining two coordinates
span one edge of that triangle factor, i.e., a segment.  Thus the $s'$
segments in the face can come from either the $s$ segment or the $t$
triangle factors of $\Pi(s,t)$.
For example,
a $\Pi(0,2)$ simplotope (see Figure \ref{tricrosstri}) can have
a $\Pi(2,0)$ face (i.e., a square).  This happens when one column in
both triangle factors is fixed at zero, and they both effectively
become segment factors, whose product is a $\Pi(2,0)$ simplotope.

Given a simplotope, we will need to know how many faces it has of a certain type.  Let $Q(a_1, \ldots, a_n; a_1', \ldots, a_n')$ be the number of $\Pi(a_1', \ldots, a_n')$ faces in $\Pi(a_1, \ldots, a_n)$.

\begin{theorem}\label{thm:generating-function}
Let $x_0 = 1$.  Then $Q(a_1, \ldots, a_n; a_1', \ldots, a_n')$ is equal to the coefficient on the $x_1^{a_1'} \cdots x_n^{a_n'}$ term in the generating function
$$
\prod_{i = 1}^n  \left( \sum_{k = 0}^i  {i+1 \choose k+1} x_k  \right)^{a_i}.
$$
\end{theorem}

\begin{proof}
%Our goal is to find the number of $\Pi(a_1', \ldots, a_n')$ faces.  
Every simplex factor in the simplotope is represented by a factor in the above product.  For example, a $3$-simplex would be represented by
$$
(x_3 + 4 x_2 + 6 x_1 + 4)
$$
This corresponds to the fact that a $3$-simplex consists of one $3$-simplex (tetrahedron), four triangles, six segments, and four vertices.   When choosing a $\Pi(a_1', \ldots, a_n')$-face of the simplotope $\Pi(a_1, \ldots, a_n)$, we must create a product of $a_1'$ segments, $a_2'$ triangles, $a_3'$ tetrahedrons, and so on, each one chosen from a factor of $\Pi(a_1, \ldots, a_n)$.  If we consider all possible ways of doing this, we get the number of $\Pi(a_1', \ldots, a_n')$-faces.   Notice those that this process is exactly analogous to finding how many ways produce a $x_1^{a_1'} \cdots x_n^{a_n'}$ term in the product above.  Hence, the coefficient on this term must be the same as the number of $\Pi(a_1', \ldots, a_n')$-faces.
\end{proof}

We now give a more explicit formula for $Q(s, t; s', t')$, which is the number of $\Pi(s',t')$ faces in
$\Pi(s,t)$. 
%Let $q$ denote the number of segment factors in a given $\Pi(s', t')$ face that come from
%segment factors in $\Pi(s,t)$.  This leaves $s'-q$ of the segment
%factors of the $\Pi(s', t')$ face to come from triangle factors of 
%$\Pi(s,t)$.

\begin{corollary} 
\label{qtheorem}
$$
Q(s, t; s', t') = {t \choose t'} \sum_{q = 0}^{s'} 2^{s-q} 3^{t-t'}
{s \choose q} {t-t' \choose s'-q} .
$$
\end{corollary}

\begin{proof}
By Theorem~\ref{thm:generating-function}, $Q(s, t; s', t')$ is simply the coefficient on the $x^{s'} y^{t'}$ term in the generating function $(x+2)^s (y + 3x + 3)^t$.   To find this coefficient explicitly, first we consider choosing exactly $t'$ copies of $y$ from the $t$ copies of $y$ in the product $(y + 3x + 3)^t$;  there are ${t \choose t'}$ ways of doing this.

From the remaining $t - t'$ factors of $(y + 3x + 3)$ and remaining $s$ factors of $(x + 2)$, we need to choose exactly $s'$ copies of $x$.  Suppose we choose $q$ of these from among the $s$ factors of $(x + 2)$.  There are ${s \choose q}$ ways of choosing $q$ copies of $x$ from the $(x+2)$ factors, and the rest of the $s' - q$ factors must be chosen from the $t-t'$ factors of $(y + 3x + 3)$.  There are $3^{s' - q} {t - t' \choose s' - q}$ ways of choosing these additional factors.   All other factors must contribute a constant:  a $2$ in the case of the $(x + 2)$ factor, and a $3$ in the case of a $(y + 3x + 3)$ factor.  This gives $2^{s - q}$ for the $(x + 2)$ factors and $3^{t - t' - (s' - q)}$ for the $(y + 3x + 3)$ factors.   Sum these all up over all possible $q$, and this gives the desired coefficient, and hence the desired number of faces.
\end{proof}

In the case $t = t' = 0$, 
the simplotope has no triangle factors.  
Because it consists only of segment factors, it is an $s$-cube. 
Using Corollary \ref{qtheorem}, we see that for this case
$$
Q(s, 0; s', 0) = 2^{s-s'} {s \choose s'},
$$
which is the formula for the surface $\mathbb{R}^{s'}$-volume of the unit $s$-cube.

\subsection*{Counting Exterior Faces.}
Consider a cover of $\Pi(a_1, \ldots, a_n)$, which consists of 
$d$-dimensional simplices, where $d = \sum_{i = 1}^n i \cdot a_i$.  Any such simplex $\alpha$ of class $c$ may or may not have,
in a given $\Pi(a_1', \ldots, a_n')$ face, an exterior $d'$-dimensional face of class $c'$, for $d' = \sum_{i = 1}^n a_i'$.
Over all $\Pi(a_1', \ldots, a_n')$ faces, $\alpha$ may have several exterior $d'$-dimensional class $c'$ faces.
Let $$F(a_1, \ldots, a_n, c; a_1', \ldots, a_n', c')$$ denote the maximum number of such faces over all possible $\alpha$.   We will be especially concerned with $F(s, t, c; s', t', c')$, which is the maximum number of $\Pi(s',t')$ faces of class $c'$ of a simplex embedded in a simplotope $\Pi(s,t)$.  Although we may not know $F$ explicitly we will derive a bound for $F$ later.

Let $V(a_1, \ldots, a_n)$ denote the largest possible class of a simplex in a cover of $\Pi(a_1, \ldots, a_n)$.
We can now formulate an inequality that a cover of $\Pi(a_1, \ldots, a_n)$ must satisfy.  This inequality 
will form the basis of a linear program that we will solve.

\begin{theorem}
\label{program-general}
Given a cover of a $\Pi(a_1, \ldots, a_n)$, 
let $x_c$ be the number of simplices of class $c$ in that cover.
Then for any tuple $(a_1', \ldots, a_n')$ and $d' = \sum_{i = 1}^n i \cdot a_i'$, we have the series of inequalities
$$
\sum_{c = 1}^{V(a_1, \ldots, a_n)} \frac{c \cdot x_c}{d'!}\  F(a_1, \ldots, a_n, c; a_1', \ldots, a_n', c) \  
\geq \  \frac{Q(a_1, \ldots, a_n; a_1', \ldots, a_n')}{\prod_{i = 1}^n i^{a_i'}}.
$$
Here, $a_n'$ varies from $0$ to $a_n$, $a_{n-1}'$ varies from $0$ to $a_{n-1}' + a_n - a_n'$, etc.
\end{theorem}

\begin{proof}This is a volume bound.  On the right side, we
  have the number of $\Pi(a_1', \ldots, a_n')$ faces $Q(a_1, \ldots, a_n; a_1', \ldots, a_n')$ multiplied by the volume 
  of these faces.  This volume is easy to calculate because it is a 
  product; each segment multiplies the volume by $1$, and each 
  triangle multiplies the volume by $1/2$, each tetrahedron multiplies it by $1/3$, etc.

The sum on the left side is an upper bound for all the exterior faces of simplices 
that could cover the volume of $\Pi(a_1', \ldots, a_n')$ by a collection of simplices, then every
$(d-1)$-facet of that face must also be covered by the same collection.
Therefore, on the left side, we only need
to count faces of simplices that are exterior facets, or that are
exterior facets of exterior facets, or exterior facets of exterior facets of exterior facets, etc.
As we will show in
Proposition \ref{facetsame}, an exterior facet will always have the
same class as the simplex itself.  Therefore we need only to consider
exterior faces that are the same class as the simplex itself, i.e., for which $c'=c$.  Each of
these faces has a volume of $\frac{c}{d'!}$.

This inequaltiy must be satisfied for any choice of $\Pi(a_1', \ldots, a_n')$.   Possible values of $a_n'$ are clearly from $0$ to $a_n$.  Possible values of $a_{n-1}'$ are from $0$ to $a_{n-1} + a_n - a_n'$, since $(n-1)$-simplices can come from the $a_{n-1}$ simplices of dimension $(n-1)$, or from the leftover $a_n - a_n'$ simplices of dimension $n$.   Similarly, $a_{n-2}'$ can vary from $0$ to $a_{n-2} + a_{n-1} + a_n - a_{n-1}' - a_n'$.  In general, $a_i'$ can vary from $0$ to $\sum_{j = i}^n a_j - \sum_{j = i+1}^n a_j'$.   
\end{proof}

Again, we will be most concerned with products of segments and triangles.   In this case, we have the corollary
\begin{corollary}
\label{program}
Given a cover of a $\Pi(s,t)$, 
let $x_c$ be the number of simplices of class $c$ in that cover.
Then for any $s', t'$ pair
$$
\sum_{c = 1}^{V(s,t)} \frac{c \cdot x_c}{(s'+2t')!}\  F(s, t, c; s', t', c) \  
\geq \  \frac{Q(s, t; s', t')}{2^{t'}}
$$
for $t'$ between $0$ and $t$, and $s'$ between $0$ and $s+t-t'$.
\end{corollary}

\subsection*{The Uniqueness of Shadow-Footprint Pairs}

In this subsection we will define shadows and footprints, and derive a result useful for counting exterior faces of simplices embedded in a simplotope.   

First we need a result about the intersection of exterior faces of a simplex.
%, and note that, unlike the other results in this section that apply to only simplotopes, this applies to any polytope.

\begin{proposition}
\label{exterior-footprints}
Let $\alpha$ be a non-degenerate simplex in a polytope
$P$, and let $\sigma, \tau$ be two exterior faces
of $\alpha$.  Then $\sigma \cap \tau$ is also an exterior face of $\alpha$.
\end{proposition}

\begin{proof}
The intersection $\sigma \cap \tau$ is simply the convex hull of the vertices they have in common.    Suppose $\sigma$ consists of $j_1 + 1$ vertices in a $j_1$-face of $P$, and $\tau$ consists of $j_2 + 1$ vertices in a $j_2$-face of $P$.   Suppose $\sigma \cap \tau$ consists of $j_3 + 1$ vertices in a $j_4$-face of $P$.  If we can show $j_3 = j_4$, this would force $\sigma \cap \tau$ to be an exterior face and the proof would be finished.

Since $\sigma \cap \tau$ is contained in a $j_4$-face of $P$, we see the $j_1$-face containing $\sigma$ and the $j_2$-face containing $\tau$ must then intersect in a face of dimension at least $j_4$.   Hence the $j_1$-face and $j_2$-face together lie in a face of dimension at most $j_1 + j_2 - j_4$.    Note also that $\sigma \cup \tau$ must have exactly $j_1 + j_2 - j_3 + 1$ points in general position, and these lie in that face of dimension at most $j_1 + j_2 - j_4$.  Hence $j_1 + j_2 - j_3 + 1 \leq j_1 + j_2 - j_4 + 1$, which gives $j_4 \leq j_3$.    Since there are $j_3 + 1$ vertices in a $j_4$ face for a simplex, we also have $j_4 \geq j_3$.  Hence $j_4  = j_3$ and $\sigma \cap \tau$ is an exterior face of $\alpha$.
\end{proof}

For $d = c_1 + \cdots, + c_n$, let $\alpha$ be an $d$-simplex of class $c \neq 0$ (non-degenerate) in the simplotope $\Delta^{c_1} \times \cdots \times \Delta^{c_n}$.
Suppose that $\sigma$ is an exterior face of $\alpha$ with dimension
$k$.  By Lemma~\ref{exterior}, the standard matrix representation of
$\sigma$ must have $d - k$ columns consisting entirely of
zeros. 
Choose one vertex $\mathbf{v}$ of $\sigma$.  Recall in our definition of class we defined a matrix $M_{\mathbf{v}}$, which starts with $M$ and then deletes every column where $\mathbf{v}$ contains a $1$, and also removes the row $\mathbf{v}$.  
Rearrange the rows of $M_\mathbf{v}$ by putting the vertices of $\sigma$ other than $\mathbf{v}$ in the first $k$ rows, followed by
the rest of the vertices of the simplotope.  Rearrange the columns so that the coordinates
that must be zero in $\sigma$ are in the final columns. The result is
$$
M_\mathbf{v} = \left [ 
\begin{tabular}{c|c} 
 & \\
\ \ \ \ A \ \ \ \ & \ \ \ zeros \ \ \ \\ \\ \hline \\ C & B \\ \\
 \end{tabular} 
\right ],
$$
where block $A$ are the nonzero columns of $\sigma$, and $C$ and $B$ are blocks corresponding to the other vertices.  Note that $A$ is a $k\times k$ block, and $B$ is a 
$(d-k) \times (d-k)$ block.  Also, observe 
that the rows that $A$ inhabits must be linearly
independent, or else $\sigma$ would not be a simplex.  Therefore we
can row reduce $M_\mathbf{v}$ to get $M_\mathbf{v}(\sigma)$ by adding multiples
of the rows of $A$ to zero out all the rows of $C$.  We are left with
$$
M_\mathbf{v}(\sigma) = \left [ 
\begin{tabular}{c|c} 
 & \\
 \ \ \ \ A \ \ \  \ & \ \ \ zeros \ \ \ \\ \\ \hline \\  \ \ \ zeros \ \ \ & B \\ \\
 \end{tabular} 
\right ],
$$
Note that $|\det[M_\mathbf{v}]|$ is the class of $\alpha$, so $|\det[M_\mathbf{v}(\sigma)]|$ is also the class of $\alpha$.  But $|\det[M_\mathbf{v}(\sigma)]| = |\det(A)\det(B)|$.  From this we can make a few observations:

\begin{proposition}
\label{facetsame}

The class of any exterior face of a simplex divides the class of the simplex.  Furthermore, the class of any exterior facet of a simplex equals the class of the simplex.

\end{proposition}

\begin{proof}
Let $\sigma$ be an exterior face.  The first statement follows noting the class of the simplex is the product $|\det(A)\det(B)|$, and the class of $\sigma$ is $|\det(A)|$. 
If $\sigma$ is a facet, then matrix $B$ must be a 1x1-matrix, with entry 1 or 0.  Since $B$ is non-degenerate, it must have determinant 1.  Thus the class of the simplex is just $|\det(A)|$, which is the class of $\sigma$. 
\end{proof}

Let $\sigma_\perp$ be the simplex spanned by the origin
and the last $d-k$ rows of $M_\mathbf{v}(\sigma)$.
Let $\pi_\sigma$ denote the linear projection of $\alpha$ onto $\sigma_\perp$
that takes the vertices of $\sigma$ to the origin and takes the last $d-k$ rows of $M_\mathbf{v}$ to the corresponding rows of $M_\mathbf{v}(\sigma)$.

\begin{proposition}
\label{pi121}

The projection $\pi_\sigma$ is one-to-one on vertices of $\alpha$ that are not in $\sigma$.  

\end{proposition}

\begin{proof}
The vertices of $\alpha$ that are not in $\sigma$ are represented by
the last $d-k$ rows of $M_\mathbf{v}$, the submatrices $B$ and
$C$.  The projection $\pi_\sigma$ simply takes $C$ and zeros it out completely,
so if $\pi_\sigma$ were not one-to-one, there would be two identical rows of $B$.
But $\det(B) \neq 0$ because the class $c\neq 0$.
\end{proof}

\begin{figure}[h]

\begin{center}
\includegraphics[scale=.75]{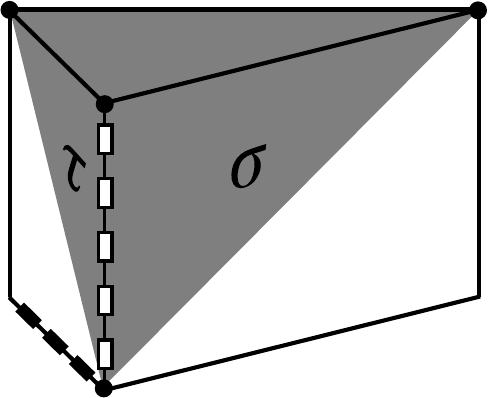}
\end{center}

\caption{\label{prism}A diagram of a shadow-footprint pair.  The simplex, shown in grey, is part of triangular prism $\Pi(1,1)$.  Both $\tau$ and $\sigma$ are faces of the simplex that are part of $\Pi(2,0)$ faces of the simplotope.  The footprint of $\tau$ with respect to $\sigma$ is shown with a white dotted line, and the shadow of $\tau$ with respect to $\sigma$ is shown with a black dotted line.}

\end{figure}

Given another exterior $k$-face $\tau$ of $\alpha$, the \emph{footprint} of $\tau$ with respect to $\sigma$ is the intersection $\sigma \cap \tau$, denoted by $\phi_\sigma(\tau)$.
The \emph{shadow} of $\tau$ with respect to $\sigma$ is $\pi_\sigma(\tau)$.
The shadow will always be a subset of $\sigma_\perp$.  The footprint will always be a subset of 
$\sigma$.

\begin{theorem}
\label{unique}
Given an exterior face $\sigma$, 
no two distinct exterior faces $\tau_1, \tau_2$ have both the same footprint and shadow with respect to $\sigma$.

\end{theorem}

\begin{proof}
Consider two exterior faces $\tau_1$ and $\tau_2$ that have the same footprint and shadow with respect to $\sigma$.  Since their footprints are the same, then $\tau_1$ and $\tau_2$ 
have the same set of vertices in common with $\sigma$.  Since their shadows are the same, by Proposition \ref{pi121}, $\tau_1$ and $\tau_2$ have the same set of vertices of $\alpha$ that are not in $\sigma$.  So $\tau_1$ and $\tau_2$ are identical.
\end{proof}

This theorem is the key insight that will allow us to bound $F(s,t,c;s',t',c')$.

\subsection*{Zero, Free, and Dependent Coordinates}

We introduce some terminology that will make subsequent proofs
clearer.  Over any subset $A$ of the simplotope, 
some of the coordinates of points in $A$ in the standard
representation may vary, and some may be fixed at $0$ or $1$ (e.g., if
$A$ is a subset of a face of the simplotope).  
Call the coordinates that are fixed at $0$ (over all of $A$)
the \emph{zero coordinates} of $A$, 
the coordinates fixed at $1$ the \emph{dependent coordinates} of $A$,
and remaining coordinates that vary over the points of $A$ the
\emph{free coordinates} of $A$.

Notice that choosing a subset of the coordinates to be zero coordinates
corresponds to a choice of face of the simplotope.

As an example, let $A$ be a $\Pi(2,0)$ face of a $\Pi(0,2)$ simplotope 
(which is a square face of the product of two triangles).
This must have one zero coordinate and two free
coordinates in each triangle factor of $\Pi(0,2)$.
Also, because there are 3 ways of picking the zero coordinate in
each factor, the 9 resulting ways of picking both zero coordinates
correspond to the 9 square faces of the product of two triangles.

For a given $\Pi(s',t')$ face, every triplet of free coordinates 
in a triangle factor of the simplotope 
contributes a triangle factor to that face.  
Likewise, every pair
of free coordinates from a segment factor of the simplotope contributes a
segment factor to the face.  Finally, a pair of free coordinates
matched with a zero coordinate from a triangle factor of the
simplotope contributes a segment factor to the face.

Now consider a non-degenerate simplex $\alpha$ in a $\Pi(s,t)$, and let
$\sigma$ be an exterior face of $\alpha$ that lives in a $\Pi(s',t')$
face of the simplotope.
Notice that the zero and dependent coordinates of a 
$\sigma$ are free coordinates of $\sigma_\perp$,
and likewise zero and dependent coordinates of $\sigma_\perp$ are free
coordinates of $\sigma$.  Thus, a segment factor of the simplotope
must correspond to either a segment factor of $\sigma$ or a segment
factor of $\sigma_\perp$.  A triangle factor of the simplotope can
correspond either to a triangle factor of $\sigma$, a triangle factor
of $\sigma_\perp$, or to a pair of segment factors: one in $\sigma$,
and one in $\sigma_\perp$.  See Figure \ref{fsdiagram}.

When taking the footprint or shadow of an exterior face $\tau$ of
$\alpha$, every free coordinate in $\tau$ corresponds to 
exactly one free coordinate in either the footprint or the shadow.

\subsection*{Combinatorial Upper Bound on $F(s, t, c; s', t', c')$}
We can determine some values of $F$.  
First, note that $F(s, t, c;s', t', c') = 0$ if $s' + 2t' > s + 2t$, 
since clearly we can't have a face with a higher dimension than that of the simplex.  
If $s = s'$ and $t = t'$, then we just have the simplex itself so $F(s,t,c;s,t,c) = 1$.  We also define all values of $F$ where $s < 0$, $t < 0$,
$s'<0$, $t'<0$, $c<1$, or $c'<1$ to be zero.  If $c > V(s,t)$ then $F$
is zero.  If $c'$ does not divide $c$, by Proposition \ref{facetsame}, $F$ is zero.

To obtain an upper bound for $F$ we need to establish some facts about parallel and tri-positioned faces.
We say that two faces of a simplotope are {\em parallel} if they have
exactly the same free coordinates.  For example, in $\Delta^2 \times \Delta^1$, the 
line connecting the points
\begin{eqnarray*}
& & (0, 1; 0, 1; 0, 0, 1)  \\
& & (0, 1; 0, 1; 1, 0, 0)
\end{eqnarray*}
is parallel to the line connecting the points
\begin{eqnarray*}
& & (0, 1; 1, 0; 0, 0, 1) \\
& & (0, 1; 1, 0; 1, 0, 0).
\end{eqnarray*}
This is because in both cases each line has two free coordinates, the 5th and 7th coordinates.  All other coordinates are either zero or dependent in both lines.

\begin{proposition}
\label{one-vertex}

Let $\alpha$ be a non-degenerate simplex in a simplotope,
and let $\sigma$ be an exterior face of $\alpha$.  Then any face of
the simplotope parallel and distinct to the one spanned by $\sigma$ can
contain at most $1$ vertex of $\alpha$.
\end{proposition}

\begin{proof}
Let $E$ be the set of points in $\alpha$ on a given face parallel
to the one spanned by $\sigma$; then $E$ has exactly the same free coordinates as $\sigma$.
Under the map $\pi_\sigma$, all the free coordinates of $\sigma$ are sent to zero, 
so therefore all the free coordinates of $E$ are sent to zero.  
Since $E$ under the map $\pi_\sigma$ no longer has any free coordinates, it is completely fixed.  So $\pi_\sigma$ maps
$E$ onto a single point.  By Proposition~\ref{pi121}, we know
$\pi_\sigma$ is one-to-one on the vertices of $\alpha$ that are not
in $\sigma$, so $E$ can contain at most one vertex of $\alpha$.
\end{proof}

This yields an immediate corollary about parallel faces.

\begin{proposition}
\label{no-parallel}

Let $\alpha$ be a non-degenerate simplex in a simplotope.
No two distinct exterior faces of $\alpha$ of dimension greater than or equal to 1
can span parallel faces of the simplotope.
\end{proposition}

We now work towards a similar result with what we will call $n$-positioned faces.

%FS: 
% use of the word span?
%
% check this and next proposition
%

\begin{lemma}
\label{lemma:n-positioned-specific}
Let $\alpha$ be a non-degenerate simplex in a simplotope, and let $S$ be a set of faces of the simplotope each spanned by vertices of $\alpha$, with $|S|=n$.  Then all three of the following cannot simultaneously hold:
\begin{enumerate}
\item the faces in $S$ have the same zero, dependent, and free coordinates outside a special factor with $n$ coordinates,
\item all $n$ of the coordinates in the special factor are free except for one zero coordinate and it is a different zero coordinate for each face in $S$, and
\item all the faces in $S$ have dimension at least $n-1$.
\end{enumerate}
\end{lemma}

\begin{proof}
Suppose there was such a set $S$.  Then outside a factor with $n$ coordinates, the faces of $S$ all have the same zero, free and dependent coordinates, and within that factor of $n$ coordinates they all have a distinct zero coordinate with all others free.  Notice this forces each face in $S$ to have the same dimension, and it is specified that this dimension is at least $n-1$.

To have a specific example in mind, consider the simplotope $\Delta^1 \times \Delta^1 \times \Delta^2$, and suppose $S$ consists of three faces (where $0$ represents a zero coordinate, $1$  a dependent coordinate, and $F$ and free coordinate):
$$
\begin{array}{rllcllclllr}
\sigma_1 =  (& F & F & ; & 0 & 1 & ; & 0 & F & F & ) \\
\sigma_2 =  (& F & F & ; & 0 & 1 & ; & F & 0 & F & ) \\
\sigma_3 =  (& F & F & ; & 0 & 1 & ; & F & F & 0 & ) \\
\end{array}
$$
In this example, $n = 3$, the special factor is the final $\Delta^2$ factor, and outside this factor the faces of the simplotope are identical in terms of the types of coordinates that define the face.

Returning to the proof, suppose that each face in $S$ has exactly $k$ vertices of $\alpha$, which means the dimension of each face in $S$ is $k-1$.  Since the dimension of each face is at least $n-1$, we have $k-1 \geq n-1$, and hence $k > n-1$. 

Let $T$ be the union of all the vertices of $\alpha$ contained in all the faces in $S$.  Since there are $n$ faces in $S$ and $k$ vertices of $\alpha$ in each face, a naive count of $T$ is $kn$ total vertices, but these could count some vertices multiple times.  Notice that no vertex can be in every face in $S$, since that would force every coordinate within the special factor  to be zero, which is impossible.  Thus each vertex is over-counted by a factor of at most $n-1$, so the number of unique vertices in $T$ is at least $\frac{kn}{n-1} = k + \frac{k}{n-1}$.   Using $k > n-1$, we have that there are at least $k+2$ unique vertices in $T$.

The dimension of a face of the simplotope is simply the number of free coordinates minus one for each factor.   The number of free coordinates in one face of $S$ is only one less than the number of free coordinates among all the faces in $S$.  Hence the dimension of the vertices in $T$ is at most $k$.   This contradicts the fact that $T$ contains at least $k+2$ unique vertices of a simplex.
\end{proof}

For $n \geq 3$ and a simplotope, we say that $\sigma_1, \ldots, \sigma_n$ faces are $n$-positioned if they have the same free coordinates except for exactly one special factor of $n$ coordinates.  Within this special factor, each face has exactly one zero coordinate, and it is a different coordinate for each face.

We can now strengthen Lemma~\ref{lemma:n-positioned-specific} by weakening requirement (1) to simply that the have the same free coordinates outside the special factor.  

\begin{proposition}
Let $\alpha$ be a non-degenerate simplex in a simplotope, and let $S$ be a set of $n$-positioned faces of dimension at least $n-1$.   Then it cannot be the case that $\alpha$ spans all the faces in $S$.
\end{proposition}

\begin{proof}
The faces in $S$ all have the same free coordinates outside of a special factor with $n$ coordinates.   Group them into $S_1, S_2, S_3, \ldots, S_t$ such that the members of each group also have the same zero and dependent coordinates outside the special factor. Notice that since different groups have different zero and dependent coordinates, no vertex in a face in $S_i$ could also be contained in a face in $S_j$ for distinct $i$ and $j$.

To have a specific example in mind, consider the simplotope $\Delta^1 \times \Delta^1 \times \Delta^1 \times \Delta^3$, and suppose $S$ consists of four faces (where $0$ represents a zero coordinate, $1$  a dependent coordinate, and $F$ and free coordinate):
$$
\begin{array}{rllcllcllcllllr}
\sigma_1 =  (& F & F & ; & 0 & 1 & ; & 0 & 1 & ; & 0 & F & F & F & ) \\
\sigma_2 =  (& F & F & ; & 0 & 1 & ; & 0 & 1 & ; & F & 0 & F & F & ) \\
\sigma_3 =  (& F & F & ; & 1 & 0 & ; & 1 & 0 & ; & F & F & 0 & F & ) \\
\sigma_4 =  (& F & F & ; & 1 & 0 & ; & 1 & 0 & ; & F & F & F & 0 & ) \\
\end{array}
$$
In this example, $n = 4$, the special factor is the final $\Delta^3$ factor, and outside this factor the faces of the simplotope form two groups:  $S_1 = \{\sigma_1, \sigma_2\}$ and $S_2 = \{\sigma_3, \sigma_4\}$, where in each $S_i$ the coordinates are identical outside the special factor.

Returning to the proof, we know by dimension requirements that $\alpha$ must have the same number of vertices in each face.  Suppose every face of $S$ has $k$ vertices from $\alpha$, and thus $k-1$ is the dimension of every face in $S$.  Since $n \geq 3$ and $k-1 \geq n-1$, we have $k \geq 3$.  Let $T$ be the union of all the vertices in the faces of $\alpha$ that span faces in $S$.  Notice $T$ must contain at least $k$ distinct vertices for each of the $t$ sets $S_i$, and so $|T| \geq tk$.

Consider the dimension of the points in $T$.  Each face starts in dimension $k-1$.  If you consider the faces within $S_i$, these combine to have dimension $k$, since only one extra free coordinate is added (within the special factor).   

I also claim that adding the points from $S_2$ to $S_1$ increases the dimension at most $1$.  We know the faces in $S_2$ has the same free coordinates as the faces in $S_1$ by assumption, so no new free coordinates are added.  The difference then is in the zero and dependent coordinates.  Since the zero and dependent coordinates are all the same within $S_2$, this means that after adding a single vertex to the affine combination, no other vertex in a face in $S_2$ will increase the size of the space.   Since only one vertex is needed to include $S_2$, this increases the dimension by at most $1$.  Similar logic holds for adding points from $S_3$, and so on for the rest of the $S_i$.  Hence, the total dimension spanned by such vertices is at most $k + t - 1$.  

Thus we have strictly more than $tk$ points on a face of dimension $k + t - 1$.  If $t = 1$, then  the result follows from Lemma~\ref{lemma:n-positioned-specific}.   Hence $t \geq 2$, and we've established $k \geq 3$. This implies that the number of points $tk$ is at least two bigger than the dimension $k + t - 1$.  We contradict the fact that $\alpha$ is non-degenerate.  
\end{proof}

The special case we focus on for $n$-positioned faces is $n = 3$, which we call {\em tri-positioned}.  Three tri-positioned faces would have the same free coordinates outside of a triangle factor, and within that triangle factor each of the three faces would have a unique zero coordinate and two free coordinates.
For example, the three square facets of a triangular prism $\Pi(1,1)$ are tri-positioned, as are the three edges of any triangular facet.

\begin{proposition}
\label{no-tri-positioned}
Let $\alpha$ be a non-degenerate simplex in $\Pi(s,t)$.
No three distinct exterior faces of $\alpha$ of dimension greater than
or equal to 2 can span tri-positioned faces of
$\Pi(s,t)$.
\end{proposition}

Now we can develop a 
combinatorial upper bound on $F$ in the case of segments and triangles.  

\begin{theorem}
\label{combbound}

If the dimension $s'+2t' \geq 2$, then 
$$
F(s, t, c; s', t', c') \leq {t \choose t'} \sum_{q = 0}^{\min(s, s')} {s \choose q} {t-t' \choose s'-q} 2^{s'-q} .
$$
Otherwise, $F(s, t, c; 1, 0, 1) \leq s+3t$ and $F(s,t,c;0,0,1) = s + 2t + 1$.

\end{theorem}

\begin{proof}
Consider a simplex $\alpha$ of the simplotope $\Pi(s,t)$.  We will
bound the number of exterior faces of $\alpha$ that are a part of a $\Pi(s',t')$ face of the
simplotope.

Recall that each of the $s'$ segment factors of the $\Pi(s',t')$ face
can arise from a segment factor of the simplotope, 
or a triangle factor of the simplotope, in which exactly two of the coordinates are not free.
Let $q$ denote the number of segment factors of the $\Pi(s',t')$ face that
arise from segment factors of the simplotope.  
Thus $s'-q$ segment factors arise from triangle factors of the simplotope.

From Proposition~\ref{no-parallel}, no
two distinct faces of $\alpha$ be in parallel faces of $\Pi(s,t)$.  This
means two different exterior faces of $\alpha$ must have a distinct
(but not necessarily disjoint) set of free coordinates.

Thus we just need to count the number of valid ways of picking
free coordinates for the faces of $\alpha$.  
Let us first choose the triangle factors of the face of $\alpha$, 
which amounts to setting all three coordinates 
in a triangle factor to be free coordinates.  
There are ${t \choose t'}$ ways of doing this.

Now let us pick the segment factors.  Notice $q$ can assume values 
from $0$, where all the segment factors come from triangles, to $\min(s,
s')$.  (Recall that $s'$ may be bigger than $s$.)
Then there 
are ${s \choose q}$ ways of choosing segment factors of the face from
the segment factors of the simplotope.  Finally, we need to pick $s'
- q$ segment factors of the face from the triangles factors of the
simplotope.  There are $t - t'$ triangles factors left over that can
be used to get segment factors of the face; we have ${t-t' \choose s'-q}$ ways of doing
this.  Notice each way of choosing a segment from a triangle could have up to
two possibilities corresponding to having a zero coordinate in two of
three places in a triangle factor.  We know we can not have a zero
coordinate in the last place in this triangle factor, because
then we would have three tri-positioned faces, and by
Proposition~\ref{no-tri-positioned}, this is impossible if the dimension $s'+2t'\geq 2$.  Therefore,
we get an extra factor of $2^{s'-q}$.  That gives the desired upper bound.

In the special case $(s', t') = (1, 0)$, Proposition~\ref{no-tri-positioned} does not hold, so
$\alpha$ may have 3 exterior edges contained in tri-positioned edges (of one triangular face), when 
there are three vertices of $\alpha$ whose factors have identical coordinates except for one triangle factor, 
and in that factor, the $3$ possibilities for edges are specified by a choice of a pair of vertices.
Other exterior edges of $\alpha$ may arise from having two vertices whose factors have identical coordinates 
except for one segment factor.  In view of Proposition \ref{no-parallel}, there are no other exterior edges of $\alpha$ possible, since
two vertices of $\alpha$ whose coordinates differ in more than one factor cannot be exterior because it will not have enough 
zero coordinates.  Since each triangle factor gives rise to at most $3$ exterior edges, and each segment factor gives rise to at most $1$ exterior edge, and every exterior edge is class $1$, we have $F(s, t, c, 1; 0, 1) \leq s+3t$.

For the special case $(s',t')=(0,0)$, a $(0,0)$-face is a vertex, and the number of vertices of $\alpha$ is exactly $s + 2t + 1$.
%Every simplex of dimension $0$ automatically has class $1$.
\end{proof}

A {\em corner simplex} is defined by a vertex $v$ of the simplotope and all the ``neighbor'' vertices connected to $v$ by edges of the simplotope.  It is the simplex with the maximum number of exterior faces in the sense that it achieves the bound given by Theorem~\ref{combbound}.

\begin{theorem}
\label{corner}

The bound given by Theorem~\ref{combbound} is achieved by a corner simplex.

\end{theorem}

\begin{proof}  Let $v$ be a vertex in the simplotope $\Pi(s,t)$; its standard coordinate representation consists of a $1$ in a single coordinate of every factor and $0$'s every entry otherwise.  Any neighbor vertex of $v$ will have the same coordinates as $v$ in every factor except one; in that factor, there will be two coordinates transposed, one containing $1$ and the other containing $0$.  We call that factor the {\em neighbor factor} for that neighbor vertex.  Let $\alpha$ be the corner simplex consisting of $v$ and all its neighbors; let $M$ be its standard coordinate representation with $v$ as the first row.

The reduced representation is obtained from the standard representation by removing one column from each factor--- in particular, the columns of $v$ that contain a $1$.  Note that these columns have at most one or two $0$'s in them (one $0$ if the column comes from a segment factor and two $0$'s if the column comes from a triangle factor), because there can be a most one or two neighbor vertices that differ from $v$ in a given factor.  
The upshot of this is the following remark: any $3$ rows chosen from $M$ will not have any zero coordinate removed when reduced relative to $v$.

Note that the reduced coordinate representation matrix $M_\mathbf{v}$ of $\alpha$ is very simple: the first row is now all zeroes and the other rows will contain exactly one $1$ and can be arranged like so:
$$
\left [ \begin{array}{ccccc}
0 & 0 & 0 & \cdots & 0 \\
1 & 0 & 0 & \cdots & 0 \\
0 & 1 & 0 & \cdots & 0 \\
0 & 0 & 1 & \cdots & 0 \\
\vdots & \vdots & \vdots & \ddots & \vdots \\
0 & 0 & 0 & \cdots & 1 \end{array} \right ]
$$
We now determine the number of exterior faces of $\alpha$ that could lie in a $\Pi(s',t')$ face of $\Pi(s,t)$.

Note that any set of $s'+2t'+1$ rows we choose will result in a $(s'+2t')$-dimensional simplex $\tau$.  If the first row is chosen, then
$\tau$ will be exterior, because the first row has $s+2t$ zero coordinates (standard or reduced), and 
each of the $s' + 2t'$ other rows will block a different coordinate from being a zero coordinate.  So there will be exactly $s+2t - (s'+2t')$ total zero coordinates (and not fewer) in $\alpha$; hence by Lemma \ref{exterior}, $\tau$ is an exterior $(s'+2t'$)-dimensional face of $\alpha$.

This construction is the only way the selection of $s'+2t'+1$ rows forming $\tau$ can be exterior as long as $s'+2t'+1\geq 3$,
because replacing the first row by any other row will decrease the number of zero coordinates by blocking a new coordinate in the neighbor factor
(see the above remark about $3$ rows).
But if the number of rows chosen is less than $3$, then replacing the first row by any other row will not change the number of zero coordinates if the two new rows have the same neighbor factor, since a zero coordinate will be created as well as destroyed.
This only occurs when either $(s',t')=(1,0)$ or $(s',t')=(0,0)$.

Consider first the case $s'+2t'\geq 2$.  Then $\tau$ (constructed above) has the same dimension as an exterior
$\Pi(s',t')$ face, but it may not lie in a $\Pi(s',t')$ face unless that face has $t'$ triangle factors and $s'$ segment factors.
To get $t'$ triangle factors in that face, we take them from triangle factors of the simplotope.  This consists of picking the two rows that have a $1$ in the triangle factor's two columns.  We do this $t'$ times, and we have $t$ triangles to choose from.  Thus there are ${t \choose t'}$ ways of doing this.
Now we only need to pick $s'$ more vertices, which will represent segment factors of the face.  These can come from triangle or segment factors of the simplotope.  Let $q$ be the number of vertices that come from segment factors.  Then there are ${s \choose q}$ ways to pick from segment factors, and ${t-t' \choose s'-q}$ ways of picking from triangle factors.  But for each vertex picked from triangle factors there are two ways of choosing it, since each triangle factor consists of two rows.  This produces an extra $2^{s'-q}$ term.  
Taking the product of these possibilities, we get the formula in Theorem \ref{combbound}.  

This estimate is too small in the case that $(s', t') = (1, 0)$.
But this case is easy to calculate.  Clearly any vertex paired with $v$ will result in $s + 2t - 1$ zero coordinates for the edge, so it will be exterior in a $1$-face of the simplotope.  This totals $s + 2t$ exterior edges so far.  We also have an exterior edge between every pair of vertices that have the same neighbor factor.  So this adds $t$ more exterior edges, for a total of $s+3t$, as desired.
\end{proof}

\section{Recurrence Relation}
\label{recurrence-relation}

In general, we expect $F(s, t, c; s', t', c')$ to be smaller when $c$ grows (big simplices don't have as many exterior faces).  The upper bound for $F$  above does not depend on $c$ at all, and therefore is only good when $c$ is small.  To get good bounds on $F$ for higher values of $c$ we use a recurrence relation based on footprint and shadow considerations.

Consider a simplex $\alpha$ of class $c$ in $\Pi(s,t)$, 
and suppose it has at
least one exterior $\Pi(s',t')$-face $\sigma$.  
Consider another
exterior $\Pi(s',t')$-face $\tau$.  We know from Theorem \ref{unique}
that every such face $\tau$ has a unique footprint/shadow combination
with respect to $\sigma$.
Hence, counting the number of footprint/shadow combinations will yield
an upper bound on the number of such faces $\tau$.  The footprint must be
exterior because intersections of exterior faces are exterior (or
empty) from Proposition \ref{exterior-footprints} below.  The shadow
must be exterior from the Proposition \ref{exterior-shadows} below.
Therefore, we count possible footprints, which is bounded by the
function of $F$ using parameters for $\sigma$, and count possible
shadows bounded by $F$ using the parameters for $\sigma_\perp$.

\begin{proposition}
\label{exterior-shadows}
In the simplex $\alpha$ of a simplotope $\Delta^{c_1} \times \cdots \times \Delta^{c_n}$, the
shadow of an exterior face $\tau$ with respect to an exterior face
$\sigma$ is an exterior face.
\end{proposition}

\begin{proof}
Again let the dimensions of $\alpha$, $\sigma$, and $\tau$ be $N$,
$n$, and $m$, respectively.  Let $p$ be the dimension of the
intersection of $\sigma$ and $\tau$, so that the dimension of the
shadow is $m-p$.  Finally, let $q$ the number of zero coordinates that
$\sigma$ and $\tau$ share in common.  Then $q+n$ is the number of zero
coordinates in the shadow of $\tau$, because it includes every
coordinate not fixed at zero in $\sigma$, as well all the coordinates
that are fixed at zero in both $\sigma$ and $\tau$.  From the
proof of Proposition \ref{exterior-footprints} 
we know the intersection has $N-p$ zero coordinates.
Also from above, $N-p = 2N - m - n - q$, which implies $q = N - m - n
+ p$, so $q+n = N - m + p$.  Therefore the number of zero coordinates
in the shadow is equal to $N - (m - p)$.  This implies that the shadow
is exterior by Lemma \ref{exterior}.
\end{proof}

\begin{figure}
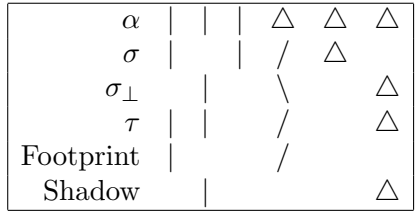

\begin{center}
\begin{tabular}{|rcccccc|} \hline
$\alpha$       & $|$ & $|$   & $|$  & $\triangle$  & $\triangle$ & $\triangle$ \\
$\sigma$       & $|$ &       & $|$   & $/$          & $\triangle$ &             \\
$\sigma_\perp$ &     & $|$ &     & $\backslash$ &             & $\triangle$ \\
$\tau$            & $|$ & $|$ &     & $/$          &                  & $\triangle$ \\   
Footprint      & $|$ &        &        & $/$          &                &             \\
Shadow         &     & $|$ &         &              &                & $\triangle$ \\ \hline      
\end{tabular}
\end{center}
\caption{\label{fsdiagram}A diagram for counting
  footprint/shadow pairs.  It illustrates how the footprint and shadow are a
  decomposition of $\tau$ into a $\sigma$ part and a $\sigma_\perp$
  part.}
\end{figure}

A good way to grasp counting footprint-shadow pairs is 
to use a footprint-shadow
diagram, like the one in Figure~\ref{fsdiagram}.  The lines, both vertical and 
slanted ones, represent segment factors, while the triangles represent triangle
factors.  In Figure~\ref{fsdiagram}, $\alpha$ is a simplex that is part of a 
cover of $\Pi(3,3)$, so it is represented by three vertical lines and three
triangles.  For the purposes of this example,
we have chosen a face $\sigma$ of $\alpha$ in a $\Pi(3,1)$ face, so 
it is represented in the diagram by three lines and one triangle;
moreover, we have chosen $\sigma$ so that one of its segment factors
arises from a triangle factor of $\alpha$, and this is represented by
the slanted line in $\sigma$'s row.
Using the diagram, it 
is easily deduced that $\sigma_\perp$ consists of two segment factors and one 
triangle factor.  Note how $\sigma_\perp$ and $\sigma$ complement each other 
with respect to $\alpha$; where there is a factor in one, 
there is not a factor in the other.  
The one exception is under the first triangle factor of $\alpha$.
Both $\sigma$ and $\sigma_\perp$ have segment factor there, but they are 
different segment factors, represented by different slants of the line.  The 
triangle factor in this case is decomposed into two segment factors.  
Using a footprint-shadow diagram can help give intuition on how to derive the following recurrence relation.

Fix $\sigma$.  Let $e$ be the number of segment factors of $\sigma_\perp$ that correspond to segment factors of $\Pi(s,t)$.
For instance, in the example of Figure \ref{fsdiagram}, $e$ is 1.  (We further explore this example near the end of this section.)
Any other simplex $\tau$ has a footprint in $\sigma$ and a shadow in 
$\sigma_\perp$, and this footprint-shadow pair is unique for each $\tau$ so we can count possible footprint-shadow pairs to get 
an upper bound on how many such $\tau$ there are.  But the footprints and shadows are exterior faces of $\sigma$ and $\sigma_\perp$ respectively,
so they may be counted recursively using $F$.  By considering all possible pairs, and summing over the potential faces the footprint can lie in, we obtain this recurrence relation:

\begin{figure}
\begin{center}
\begin{tabular}{|rcccccc|} \hline
$\alpha$       & $|$ & $|$   & $|$  & $\triangle$  & $\triangle$ & $\triangle$ \\
$\sigma$       & $|$ &       & $|$   & $/$          & $\triangle$ &             \\
$\sigma_\perp$ &     & $|$ &     & $\backslash$ &             & $\triangle$ \\
$\tau$            & $|$ & $|$ &     & $\triangle$          &                  & $/$ \\   
Footprint      & $|$ &        &        & $/$          &                &             \\
Shadow         &     & $|$ &         &  $\backslash$      &            & $/$ \\ \hline      
\end{tabular}
\end{center}
\caption{\label{fsdiagram2} Same $\alpha$ and $\sigma$ as in Figure \ref{fsdiagram} but a different $\tau$.
Note that this $\tau$ has the same footprint as the $\tau$ in Figure \ref{fsdiagram}, but the shadow of $\tau$ is different.
}
\end{figure}

\begin{theorem}
\label{recurrence}
The quantity $F(s, t, c; s', t', c')$ is less than or equal to
$$
%\max_{e}   \sum_{w=0}^{w'(e)} \sum_{k=1}^{c'}  \sum_{j=0}^{t'} \sum_{i=w}^{i'(j,w)} 
\max_{e}   \sum_{w=0}^{w'} \sum_{k=1}^{c'}  \sum_{j=0}^{t'} \sum_{i=w}^{i'} 
F\left(s',t',c',i,j,k \right) 
\cdot
% F\left(s''(e),t''(e),\frac{c}{c'},s'-i+2w,t'-j-w,\frac{c'}{k}\right) .
 F\left(s'',t'',\frac{c}{c'},s'-i+2w,t'-j-w,\frac{c'}{k}\right) .
$$
%where $s''(e)= s'-s+2e$, $t''(e)=s+t-s'-t'-e$,  $w'(e)=\min(s'-s+e,t')$, $i'(j,w)=\min(s'+t'-j,s'+w)$,
where $s'', t''$ and $w'$ are functions of $e$ and $i'$ is a function of $j$ and $w$ via these relations:
$s''= s'-s+2e$, $t''=s+t-s'-t'-e$,  $w'=\min(s'-s+e,t')$, $i'=\min(s'+t'-j,s'+w)$.
The maximum is taken over $e$ ranging from $\max(0, s-s')$ to $\min(s+t-s'-t',s)$.  

For purposes of this recursion, if $s + 2t \geq 2$, then we will take $F(s, t, c; 1, 0, 1) \leq s + 2t$ and $F(s, t, c; 0, 0, 1) = 1$.
\end{theorem}

There's a lot going on in this recursion which will be explained as we proceed through the proof. 

\begin{proof}
Fix a simplex $\alpha$ in $\Pi(s,t)$, and suppose $\alpha$ is of class $c$.
Recall that $F(s,t,c;s',t',c')$ counts the number of faces of $\alpha$ in a certain family, namely, 
the number of faces of $\alpha$ that live in a $\Pi(s',t')$ face of $\Pi(s,t)$ and that are of class $c'$.

We fix one such face of this family, call it $\sigma$.  (If there is no such face, the above inequality will trivially hold.)
Any other $\tau$ in this family will have a unique footprint-shadow pair with respect to $\sigma$, so it suffices count
the number of possible footprint-shadow pairs, 
where the footprint is an exterior face of $\sigma$ and the shadow is an exterior face of $\sigma_\perp$.

To do this, we consider footprints that lie in some $\Pi(i,j)$ face and are of class $k$,
then count the number of shadows that could be associated with this footprint, then sum up over all possible $i,j,k$.  Note that this depends on an additional quantity $e$ corresponding to the number of segment factors of $\sigma_\perp$ that correspond to segment factors of $\Pi(s,t)$.  
Since not all $e$ may be achieved by some $\sigma$, we take the maximum over all $e$ as a bound.
%Taking the minimum over all possible values of $e$ corresponds to picking the best choice of $\sigma$ with respect to which we take
%footprints and shadows.

The number of possible footprints that lie in some $\Pi(i,j)$ face and are of class $k$ is less than or equal to:
$$
F(s', t', c'; i, j, k) .
$$
%\label{footprintbound}
This follows from the definition of $F$, noting that a footprint is an exterior face of $\sigma$ in a $\Pi(i,j)$ face, and $\sigma$ is a simplex in a 
$\Pi(s',t')$ face of class $c'$.

For each such footprint, there are a number of shadows that could be associated with this footprint, and bounding this number is complicated 
by the fact that the number of segment and triangle factors of a shadow is not uniquely determined by a footprint.  See, for instance, the $\tau$ in Figures \ref{fsdiagram} and \ref{fsdiagram2}, which have the same footprint dimensions $i$ and $j$ but shadows with different numbers of segment and triangle factors. 
These numbers depend on a quantity $w$, which we define to be the number of triangle factors of $\tau$ that are chosen from 
triangle factors of $\alpha$ that support segment factors of $\sigma$.  

Recall that $e$ was defined above as the number of segment factors of $\sigma_\perp$ that correspond to segment factors of $\Pi(s,t)$.  Then
$s-e$ is the number of segment factors of $\alpha$ that support segment factors of $\sigma$, so 
the number of triangle factors of $\alpha$ that support segment factors of $\sigma$ is $s'-(s-e)$.  
Thus $w$ varies between $0$ and $s'-s+e$.  Then for a given footprint, the number of shadows is bounded by
%\label{shadowbound}
$$
 \sum_{w=0}^{s'-s+e}
 F\left(s'',t'',\frac{c}{c'},s'-i+2w,t'-j-w,\frac{c'}{k}\right).
$$
This bound follows from the definition of $F$, noting that a shadow of $\tau$ is an exterior face of $\sigma_\perp$, 
which lives in a $\Pi(s'',t'')$ simplotope face $H$ of $\Pi(s,t)$. 
What remains is to figure out is what are the appropriate dimensions of $\sigma_\perp$ and the shadow, 
and what kind of simplotope faces they live in.

From the argument above, we know the number of triangle factors of $\alpha$ that support segment factors of $\sigma$ (hence also $\sigma_\perp$)
is $s'-s+e$, the total number of segment factors of $H$ is obtained by adding $e$ to this, hence $s''=s'-s+2e$.
The number of triangle factors of $H$ is $t$ (the number of triangle factors of $\alpha$) minus $t'$ (the number of triangle factors used by $\sigma$) minus $s'-s+e$ (the number of triangle factors of $\sigma$ supporting segment factors of $\sigma_\perp$).
Therefore $t''=t-t'-(s'-s+e) = s+t-s'-t'-e$.

The class of $\sigma_\perp$ is the class of
$\alpha$ divided by the class of $\sigma$, i.e., $c/c'$.

The total dimensions of the shadow and footprint must sum to the dimension of $\tau$, which lives in a $\Pi(s',t')$ face, and if this were true individually for the segment factors and triangle factors by themselves, then the number of segment and triangle factors of the shadow would be $s'-i$ and $t'-j$, respectively.  
However, as noted above, in some instances a segment supporting the footprint face (a subface of $\sigma$) may combine with a segment supporting the shadow face (a subface of $\sigma_\perp$) to form a triangle factor of $\alpha$.   This can only happen when $\sigma$ and $\sigma_\perp$ are on supporting faces that have segment factors from the same triangle factor of $\alpha$; we saw there are at most $s'-s+e$ such triangle factors.  
Note that $w$, as defined above, is the number of factors of $\alpha$ for which the footprint and shadow have segment factors from the same triangle factor of $\alpha$.  Then $w$ can vary between $0$ (when $\tau=\sigma$) and $\min(s'-s+e,t')$, and for a given $w$, there are $2w$ more segments in the footprint and shadow dimensions, and $w$ fewer triangles.  Thus, the number of segment factors of the
shadow must be $s'-i+2w$, so that when we combine this with the $i$
segments of the footprint and the remove of the segments caused by splitting factors of $\alpha$, we have $s'$ total segment factors.
Likewise, the shadow must have $t'-j-w$ triangle factors.

Finally, the
class of the shadow and footprint must multiply to the class of
$\tau$, so the class of the shadow must be $c'/k$.

Now we sum over all possible $i$, $j$, and $k$.
As long as we count every possible footprint-shadow combination,
we will have an upper bound on the number of exterior faces $F$.

The index $j$ can run from $0$ to at most $t'$, the number of triangle factors supporting $\sigma$.
The index $i$ is at least $w$ because there are at least that many segment factors supporting $\sigma$.  The index $i$ cannot be more than 
$s'+t'-j$, the number of segment and triangle factors supporting $\sigma$ minus the number of triangles used by $j$, and $i$ cannot also be more than $s'+w$, the number of segments supporting $\alpha$ together with the number of triangle factors of $\alpha$ that support segment factors of $\sigma$.
Hence $i$ is at most $\min(s'+t'-j,s'+w)$.
The index $k$ can run from $0$ to $c'$ (in fact $k$ must be a divisor of $c'$ but that is already accounted for in the initial conditions for $F$).

What are the possible values of $e$?  Recall that $e$ is the number of segment factors of $\Pi(s,t)$ that
are segment factors of $\sigma_\perp$.  The smallest $e$ occurs in the
case where all $s'$ segment factors of $\sigma$ are from segment
factors of $\Pi(s,t)$; this leaves either $0$ or $s-s'$ segment factors for
$\sigma_\perp$, whichever is larger.  The largest $e$ occurs in the case where as many as possible of
the $s'$ segment factors of $\sigma$ come from triangle factors of
$\Pi(s,t)$.  The number of triangle factors available is at most $t-t'$, so this leaves at least $\max(s'-t+t',0)$ segments of $\sigma$ to come from segment factors of $\Pi(s,t)$, or at most $s-\max(s'-t+t',0)=\min(s-s'+t-t',s)$ segment factors from which $\sigma_\perp$ may be supported. 

We must justify why, if $s + 2t \geq 2$, we can use $F(s, t, c; 0, 0, 1) = 1$ and $F(s, t, c; 1, 0, 1) \leq s + 2t$ for the purposes of the recursion.  If, in the recursion, we encounter $F(s, t, c; 0, 0, 1)$ for counting footprints, then we are asking: how many different exterior faces, with the same shadow, could have a  single vertex footprint?  The answer is at most one, since two exterior faces with the same shadow but different single vertex footprints would be in parallel faces.  And if we encounter $F(s, t, c; 0, 0, 1)$ for counting shadows, we are asking: how many different exterior faces, with the same footprint, could have a shadow as a single vertex?  Once again, if there were two, they would have to be in parallel faces of the simplotope. 

Similarly, if we encounter $F(s, t, c; 1, 0, 1)$ in the recursion for counting footprints, then how many different exterior faces, with the same shadow, could have an edge as a footprint?  Well, the answer is $s + 2t$, not $s + 3t$, since if all three edges could be the footprint in a triangle factor, that would indicate three tri-positioned faces, which we know cannot all be exterior.   And if we encounter $F(s, t, c; 1, 0, 1)$ to count shadows, only two of the three possible edges of a triangle factor can be used without yielding tri-positioned faces.
\end{proof}

Notice that when $t=t'=0$, then $e$ is fixed at $s-s'$ and $j, w$ are fixed at zero, so we have
$$
F(s, 0, c; s', 0, c') \leq \sum_{i, j, k} F(s', 0, c'; i, 0, k) \cdot F(s-s', 0, \frac{c}{c'}; s'-i, 0, \frac{c'}{k} ),
$$
the same recurrence relation derived by Bliss and Su for cubes~\cite{blsu05}.

\subsection*{Examples}

To understand Theorem \ref{recurrence}, it is helpful to consider examples.

First, consider a prism $\Delta^1 \times \Delta^2 = \Pi(1,1)$.  Suppose we want an upper bound on the number of exterior faces
a class 1 simplex $\alpha$ can have inside a square face of the prism. 
Then $s=1, t=1, s'=2, t'=0$ and $c = c' =1$ and we want to compute $F(1, 1, 1; 2, 0, 1)$.  The first three
numbers indicate that we have a class 1 simplex $\alpha$ in 
$\Pi(1,1)$, 
and the last three numbers indicate that we are looking for exterior
faces of $\alpha$ in $\Pi(2,0)$ faces of the prism (square faces, see Figure~\ref{prism}).  From
Theorem~\ref{combbound}, we know that the maximum number of exterior
faces is $2$,  shown in Figure~\ref{prism} as $\sigma$ and $\tau$, square faces of a corner simplex $\alpha$.  
As we will see, our recurrence relation produces an upper bound of $3$.

By making appropriate substitutions, we see in Theorem~\ref{recurrence} that the following indices are fixed: $j=0$, $k=1$, $w=0$, and $e=0$ (because $\sigma$ lies on a square face and must use the segment factor of $\alpha$, so $\sigma_\perp$ cannot).  Then we have
\begin{eqnarray*}
F(1, 1, 1; 2, 0, 1) &\leq&
\sum_{i=0}^{2} F(2,0,1;i,0,1) \cdot F(1,0,1;2-i,0,1) \\
& \leq & F(2, 0, 1; 0, 0, 1) F(1, 0, 1; 2, 0, 1) \\
& + & F(2, 0, 1; 1, 0, 1) F(1, 0, 1; 1, 0, 1) \\
& + & F(2, 0, 1; 2, 0, 1) F(1, 0, 1; 0, 0, 1).
\end{eqnarray*}
We see that $F(1, 0, 1; 2, 0, 1)=0$ because there cannot be a square face of a line segment. 
In the second term,
$F(2, 0, 1; 1, 0, 1)=2$ since a simplex in a square can have two exterior faces that are segments.
In the third term, $F(2, 0, 1; 2, 0, 1)=1$ by definition. 
We know $F(1, 0, 1; 1, 0, 1)=1$ and $F(1,0,1;0,0,1)=1$ by Theorem \ref{combbound}.
Hence $F(1,1,1;2,0,1) \leq 0+2\cdot 1 + 1 \cdot 1 = 3$.

Let us look at how this is interpreted geometrically, using
Figure~\ref{prism}.  We count footprints and shadows with respect to
$\sigma$, a $\Pi(2,0)$ face we picked arbitrarily. The first term
above corresponds to trying to get an entire square face from the
shadow.  As in Figure~\ref{prism}, the shadow is just a segment, so it
cannot have a square face, therefore the first term is zero.  The
second term corresponds to using a segment exterior face of both the
shadow and the footprint.  There are two ways to do this, so this term
is equal to two.  The final term corresponds to ways in which the footprint is in a square face.
%This corresponds to when $\tau=\sigma$ but the shadow may be chosen in two ways?

Consider another example: the triangle cross a square $\Pi(2,1)$, 
shown in Figure~\ref{tricrosssquare} below.  We know $F(2, 1, 1; 1, 1, 1)$ 
is, for a class $1$ simplex $\alpha$ in $\Pi(2,1)$, the maximum number of class $1$ facets of $\alpha$ that live in prism faces.  
Since  $s = 2, t = 1, s' = 1, t' = 1$, then in the recurrence we have $e=1$, $w=0$, $k=1$ and $j$ runs from $0$ to $1$, and $i$ runs from $0$ to $1$.
So there will be four terms in the sum:
\begin{eqnarray*}
F(2, 1, 1; 1, 1, 1) 
&\leq& 
 \sum_{j=0}^{1} \sum_{i=0}^{2-j} 
F(1,1,1;i,j,1) \cdot F(1,0,1;1-i,1-j,1)\\
&\leq
      & F(1, 1, 1; 0, 0, 1)F(1, 0, 1; 1, 1, 1) \\
& + & F(1, 1, 1; 1, 0, 1)F(1, 0, 1; 0, 1, 1) \\
& + & F(1, 1, 1; 0, 1, 1)F(1, 0, 1; 1, 0, 1) \\
& + & F(1, 1, 1; 1, 1, 1)F(1, 0, 1; 0, 0, 1)\\
&= & 1 \cdot 0 + 4 \cdot 0 + 1 \cdot 1 + 1 \cdot 1 = 2.
\end{eqnarray*}

Finally, for a larger example, examine one term in the recursion
for $F(3, 3, 2; 3, 1, 1)$.  Here, $\alpha$ is a class $1$ simplex in a
cover of 
$\Pi(3,3)$, and we are searching for exterior $\Pi(3,1)$ faces of class $1$.  This example is illustrated in Figure~\ref{fsdiagram}.  With the choice of $\sigma$ in the figure, we see that $e = 1$ since $\sigma_\perp$ has one segment factor that comes from a segment factor in $\alpha$.  
The face $\tau$ is counted in the product $F(3, 1, 1; 2, 0, 1) F(2, 1, 1; 1, 1, 1)$, when $s'=3$, $t'=1$, $i=2$, $j=0$, and $w=0$.
We see the triple $(3, 1, 1)$  represents $\sigma$, the triple 
$(2, 0, 1)$  represents the footprint, the $(2, 1, 1)$ represents $\sigma_\perp$, and the $(1, 1, 1)$ represents the shadow.
Figure \ref{fsdiagram2} represents the same situation but with a different $\tau$; this $\tau$ has $w=1$ and is counted in the product
$F(3, 1, 1; 2, 0, 1) F(2, 1, 1; 3, 0, 1)$.
%, when $s'=3$, $t'=1$, $i=2$, $j=0$, and $w=1$.

\subsection*{A Linear Program}

Using the inequalities derived in Theorem \ref{program}, we can use a linear program
to minimize the number of simplices that meet the
requirements of the inequalities, i.e.,:
$$
\min \sum_{c=1}^{c=V(s,t)} x_c  \text{ subject to the inequalities in Theorem \ref{program}.}
$$

To get the best results, we need good bounds on $V(s, t)$, which is
the largest possible class of a simplex in a cover of $\Pi(s,t)$.
This is a hard problem, related to the Hadamard maximum determinant
problem.  
We use known values for cubes of small dimension as upper bounds on the values for
simplotopes, because $\Pi(s,t)$ can be seen as a subset of a
$s+2t$-cube, so the largest simplex in $\Pi(s,t)$ must be less than
or equal to the largest simplex in a $s+2t$-cube.

However, for many values of $(s, t)$, this is not a very good bound.
Therefore for specific low dimensions, we ran a brute force computer
program that checked the classes of all possible simplices that could
be part of a triangulation and returned the highest one.  
Beyond the trivial $V(1, 1) = 1$, $V(0, 2) = 1$, 
we find $V(2, 1) = 2$, $V(1, 2) = 3$, and $V(0, 3) = 4$.  
Although this improvement is only for low
dimensions of $V(s,t)$, it can affect higher dimensions of
our bounds because of the recursion that takes place in determining
$F$.

In determining values of $F$, we used upper bounds that result from
the recurrence relation, using the combinatorial upper bound on $F$ as
a base case of the recurrence or anytime it gave a better bound.

%Using \texttt{lp\_solve}, 
Using SAGE Math,
we solved the linear program to obtain the values shown in Table \ref{results-table}.
Our notebook for this calculation is available online~\cite{seacrest17}.   %at \url{https://cloud.sagemath.com/projects/4dd9c3c7-9a3b-4e4c-9989-1f44bdc7d6ca/files/Seg-Tri-Lower-Bounds.sagews}.

\begin{figure}[h]
\begin{center}
\includegraphics[scale=.5]{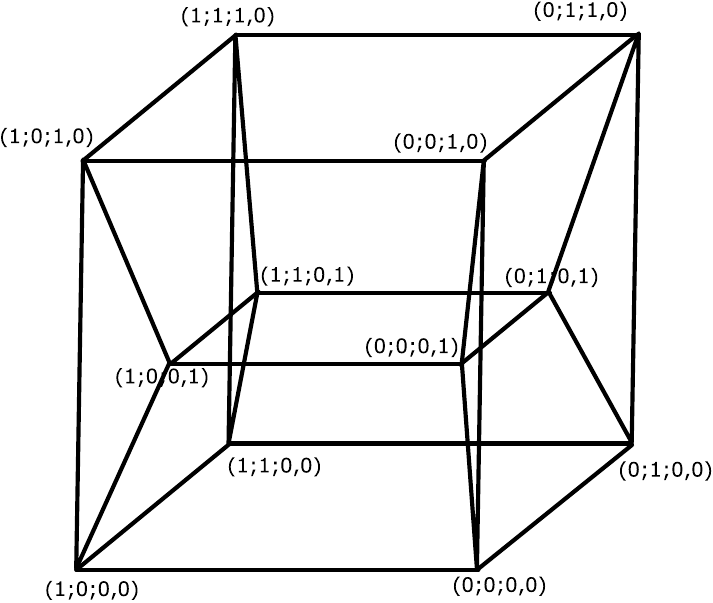}
\end{center}
\caption{\label{tricrosssquare}
A Schlegel diagram of a triangle cross a square, labeled by a reduced 
coordinate representation.}
\end{figure}

\section{The Triangle Cross Square}
\label{tri-cross-sq}

In this section, we determine a minimal triangulation of $\Pi(2,1)=\Delta^1 \times \Delta^1 \times \Delta^2$, the product of a square and a triangle.
In this special case, we can check our program bounds and, in fact,
improve them.  Our lower bound from the linear program 
shows that at least $9$ simplices are required to cover  $\Pi(2,1)$.  However, we
shall show that no cover of $9$ simplices exists, using 
the following two key propositions:

\begin{proposition}
\label{largest-is-two}
Every simplex in a cover of $\Pi(2,1)$ has an exterior facet, and
therefore the largest class of a simplex in a cover of $\Pi(2,1)$ is
class $2$.
\end{proposition}

\begin{proof}
Consider how 5 vertices of $\Pi(2,1)$ can be chosen as corners of a
simplex.  For each, the triangle factor must be one of 
$(0, 0, 1)$, $(0, 1, 0)$, or $(1, 0, 0)$ in standard coordinates.
Each of these three choices must be picked at least once
to get 5 affinely independent points, so there must be one of these choices that is chosen at
most once, hence the other two are chosen exactly 4 times total. 
The corresponding four points have a zero coordinate in common, which means, by Theorem
\ref{exterior}, the simplex has to have an exterior facet.  
By Theorem \ref{facetsame}, the class of the simplex must be equal to the class of this facet.
The exterior faces of $\Pi(2,1)$ are $\Pi(1,1)$ and $\Pi(3,0)$ simplotopes.
The largest class that a simplex can have in either of these is $2$,
which occurs for a $\Pi(3,0)$ face (and cannot occur for a
$\Pi(1,1)$ face).
\end{proof}

The \emph{center} of $\Pi(2,1)$ is the point fixed by all isomorphisms of $\Pi(2,1)$ with itself;
 in reduced coordinates is $(1/2; 1/2; 1/3, 1/3)$.

\begin{proposition}
\label{contains-center}
Any class $2$ simplex of a cover of $\Pi(2,1)$ will contain the
center of $\Pi(2,1)$ in the interior of a facet of the simplex.
\end{proposition}

\begin{proof}
From the proof of Proposition~\ref{largest-is-two}, we know that a class $2$
simplex must have an exterior facet of class $2$ in one of the
$\Pi(3,0)$ facets of the simplotope.  Therefore, every class $2$
simplex must consist of one of these class $2$ facets coned to a
point not on that facet.  There are 24 possible class $2$ simplices,
because there are two ways of choosing a class $2$ simplex facet from
a $\Pi(3,0)$ facet, 
there are three possible $\Pi(3,0)$ facets to choose from,
and there are four ways of picking a point not on that facet once it
is chosen.  However, we now show that each of these $24$ class $2$ simplices
is isomorphic up to rotations and reflections about the center.

Fix a class $2$ simplex and consider its standard matrix
representation $M$.  Permuting the columns associated to each factor 
corresponds to rotations and reflections of the simplotope that leave
the center fixed.  Furthermore, every permutation 
corresponds to a different class $2$ simplex 
because otherwise two columns of $M$ would be identical, which cannot
occur if the simplex is non-degenerate.
There are $24= 2! \times 2! \times 3!$ such permutations that permute columns of $M$ within each factor.
These account for all of the $24$ class $2$ simplices determined above.

Since the class $2$ simplices are isomorphic, without loss of generality, we can check that the proposition holds for one class $2$ simplex.
Consider the following class $2$ simplex spanned by these vertices in reduced
coordinates:  $(0; 0; 0, 0)$, $(0; 1; 0, 1)$, $(1; 0; 0, 1)$, $(1; 1; 0, 0)$, $(0; 0; 1, 0)$.
%\begin{eqnarray*}
%&(0; 0; 0, 0)& \\
%&(0; 1; 0, 1)& \\
%&(1; 0; 0, 1)& \\
%&(1; 1; 0, 0)& \\
%&(0; 0; 1, 0)&
%\end{eqnarray*}
Suppose some convex combination of these points produced the center point using coefficients $a_1$, $a_2$, $a_3$, $a_4$, and
$a_5$ respectively.  Such a combination solves the system of equations
$$
\left [ \begin{array}{ccccc} a_1 & a_2 & a_3 & a_4 & a_5 \end{array}
\right ]  \left [ \begin{array}{cccc} 0 & 0 & 0 & 0 \\ 0 & 1 & 0 & 1
    \\ 1 & 0 & 0 & 1 \\ 1 & 1 & 0 & 0 \\ 0 & 0 & 1 & 0 \end{array}
\right ]  = \left [ \begin{array}{cccc} 1/2 & 1/2 & 1/3 & 1/3
  \end{array} \right ].
$$

Furthermore, as a convex combination we know $a_1 + a_2 + a_3 + a_4 + a_5 = 1$.  
%We arrive at the system of linear equations:
%\begin{eqnarray*}
%a_1 + a_2 + a_3 + a_4 + a_5 & = & 1. \\
%a_3 + a_4 & = & 1/2. \\
%a_2 + a_4 & = & 1/2. \\
%a_5 & = & 1/3. \\
%a_2 + a_3 & = & 1/3.
%\end{eqnarray*}
These equations can be solved to find $a_1 = 0$, $a_2 = 1/6$, $a_3 = 1/6$, $a_4 = 1/3$, and $a_5 = 1/3$.  
Since exactly one of the coefficients is zero,
the center lies on a three dimensional face of the simplex
but no two dimensional face, so the center must be interior to a facet of the simplex.
\end{proof}

This yields an immediate corollary by noting that three half-spaces through one point must have an overlapping pair:

\begin{corollary}
\label{overlap}
Given any three class $2$ simplices in $\Pi(2,1)$, 
at least two of the three simplices must overlap, i.e., 
there is a point interior to both.
\end{corollary}

With these, we obtain the main result of this section.

\begin{theorem}
The minimal cover of $\Pi(2,1)$ is $10$ simplices, and this can be
accomplished with a triangulation.
\end{theorem}

\begin{proof}
Proposition \ref{largest-is-two} shows that the largest class of a simplex in $\Pi(2,1)$ is
two.   Table \ref{results-table} shows that a cover must have at least $9$ simplices.  
Then there must be at least three class
$2$ simplices in the cover; otherwise, two class $2$ simplices and
seven class $1$ simplices would not be enough to cover $\Pi(2,1)$
(which has total class 12).

On the other hand, there must be at least six class $1$ simplices in any cover, because any pair of opposite prism facets of 
$\Pi(2,1)$ require six tetrahedra to cover them that are facets of simplices in $\Pi(2,1)$.  These are distinct simplices,
because no simplex can have exterior facets in parallel facets of $\Pi(2,1)$, by Proposition \ref{no-parallel}.

Thus a size 9 cover must consist of three class 2 simplices and six class
1 simplices.  Since their total class is 12, these would have to cover
$\Pi(2,1)$ without overlapping interiors; however,
Corollary~\ref{overlap} shows that this is impossible.  Therefore, a
cover of $\Pi(2,1)$ must be size 10 or more.

%FS: reduced?
We now exhibit a size 10 triangulation, then describe its construction.  Order the vertices of $\Pi(2,1)$ numerically by their reduced coordinates, and label them by the twelve symbols $1,2,3,4,...,8,9,0,\#,*$.  Thus: $1$ represents $(0;0;0,0)$, $2$ represents $(0;0;0,1)$, $3$ represents $(0;0;1,0)$, $4$ represents $(0;1;0,0)$, and so on, until finally $*$ represents $(1;1;1,0)$.   

Using these symbols, we now specify 10 simplices of a minimal triangulation of $\Pi(2,1)$ by their $5$ vertices: 
\begin{eqnarray*}
[1850*],&1450*, 1456*, 1356*,&
[1358*], 1398*, 1798*, 1708*,\\
 \#850*,& &13582. 
\end{eqnarray*}
Two simplices are said to be \emph{adjacent} if they meet face-to-face along a facet or, equivalently, if they differ in just one vertex.
The 8 simplices in the first row above form a cycle--- each is adjacent to the simplices next to it in the displayed row, with row ends 1850* and 1708* also adjacent.  The bracketed simplices $1850*$ and $1358*$ play a special role in this triangulation--- they are the ``fat'' simplices of class 2, while all other simplices are class 1.  In addition $1850*$ and $1358*$ are adjacent to each other along the common facet $158*$, and they are each adjacent to one of the simplices in the bottom row:
$\#850*$ is the ``corner'' simplex at vertex $\#$ and is adjacent to $1850*$, 
and $13582$ is the corner simplex at vertex $2$ and is adjacent to $1358*$.  These are all the adjacencies in the triangulation. 
Every other facet of a simplex lies in an exterior facet of the simplotope.

\begin{figure}[h]
\begin{center}
\includegraphics[scale=.6]{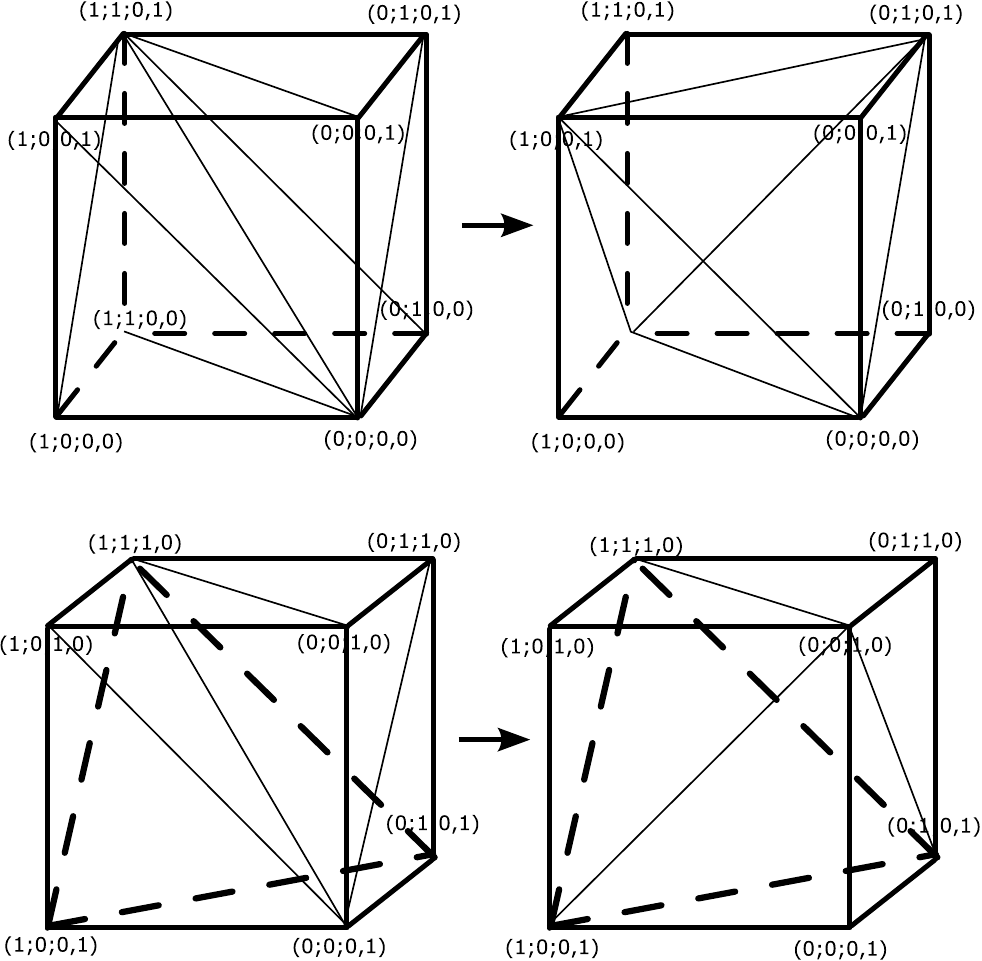}
\end{center}
\label{flips}
\caption{
Two replacements used in our construction of a $10$ simplex cover of
$\Pi(2,1)$.  
Cones over these $3$-dimensional simplices produce $4$-dimensional simplices that are part of a triangulation of $\Pi(2,1)$.}
\end{figure}

Here is how we constructed this triangulation.  Start with the standard triangulation of $\Pi(2,1)$ consisting of $12$ simplices, 
where every simplex is spanned by the two vertices $(0;0;0,1)$ and $(1;1;1,0)$ and three other vertices. 
Consider the convex hull of the following subset of points of $\Pi(2,1)$:
\begin{eqnarray*}
\ & \ & (0; 0; 0,0), (0; 1; 0, 0), (1; 0; 0, 0), (1; 1; 0, 0),\\
\ & \ &  (0; 0; 0, 1), (0; 1; 0, 1),  (1; 0; 0, 1), (1; 1; 0, 1), (1; 1; 1, 0).
\end{eqnarray*}
This is a cone over a cube, with apex at $(1; 1; 1, 0)$.    
The standard triangulation of $\Pi(2,1)$ triangulates this cone by
six simplices, which arise as a cone over a standard triangulation of the
$3$-dimensional cube, depicted at top left of Figure \ref{flips}.
The triangulation $T$ of this cone can be replaced with a new triangulation $T'$ formed by the cone of $(1;1;1,0)$ 
over a triangulation of the cube of size 5; see top right of Figure \ref{flips}.
This produces a size $11$ cover of $\Pi(2,1)$.  In fact, it is a triangulation, because both $T$ and $T'$ meet the rest of the triangulation in 
exactly the same way, i.e., the facets of $T$ and $T'$ that are on the boundary of the cone but not exterior to $\Pi(2,1)$ are exactly the same.
Referring to the top left diagram of Figure \ref{flips}, those facets are cones of $(1;1;1,0)$ over any of the triangles on the three faces of the cube containing $(0;0;0,0)$.  Note that these triangles are unchanged in the top right diagram of Figure \ref{flips}.

Now, consider the convex hull of the following points:  
\begin{eqnarray*}
\ & \ & (0; 0; 1, 0), (0; 1; 1, 0), (1; 0; 1, 0), (1; 1; 1, 0), \\
\ & \ & (0; 0; 0, 1), (0; 1; 0, 1), (1; 0; 0, 1), (0; 0; 0, 0).
\end{eqnarray*}
This is a cone of $(0;0;0,0)$ over a cube with a corner cut off; see the bottom left diagram of Figure \ref{flips}.
In the size $11$ cover of $\Pi(2,1)$ above, this 
cone has a triangulation $W$ by $5$ simplices, and it is the cone over the triangulation depicted in the bottom left diagram of Figure \ref{flips}.
However, we can replace $W$ by $W'$, the cone over the triangulation in the bottom right diagram of Figure \ref{flips}, which has just 4 simplices.   
The replacement of $W$ by $W'$ produces a triangulation, because the facets of $W$ and $W'$ that are on the boundary of the cone but not exterior to $\Pi(2,1)$ are exactly the same.  In the bottom left diagram of Figure \ref{flips}, such facets are formed by the dotted triangle, and the two vertical facets bordering the dotted triangle.  Since these facets are unchanged in the bottom right diagram, we obtain the triangulation described above
that has just $10$ simplices.

\ignore{******
standard
$$
125#* 
128#*
178#*
170#*
140#*
145#*
$$
$$
1256*
1236*
1239*
1289*
$$
$$
1456*
1789*
$$
-----------}

\end{proof}

\ignore{*********************************
\section{Improving the Bound}
\label{improving}

Corner simplices are especially important as they have the greatest
number of exterior faces in the sense of Theorem~\ref{corner}.  By
considering the fact that there can only be so many corner simplices,
we can improve the bounds given by our linear program.

\begin{theorem}
Given a simplex in $\Pi(s,t)$ with ${t \choose t'} \sum_{q = 0}^{\min(s, s')} {s \choose q} {t-t' \choose s'-q} 2^{s'-q}$ exterior $s', t'$, the simplex must be a corner simplex.
\end{theorem}

\begin{proof}
If a simplex has that many exterior vertices, then in some sense every possible combination of $s'$ segment factors and $'t$ triangle factors is an exterior face.  Therefore choosing the first $s'$ segment factors and the first $t'$ triangle factors corresponds to a face.  Likewise, choosing the $s'$ segment factors skipping the first one, and the $t'$ triangle factors skipping the first one corresponds to a face.  From the proof of Proposition~\ref{exterior-footprints}, we know that $q = N - n - m + p$.  We know these two faces holds $s'-1$ segment factors and $t'-1$ triangle factors the same.  Therefore, $p = s'+2t'-3$.  We also know $n = s'+2t'$, $m = s'+2t'$, and $N = s+2t$.  Therefore, $q = s+2t-s'-2t'-3$, which is the number of zero coordinates the two faces have in common.  
\end{proof}

********************}

\section{Discussion}
Our results in Table \ref{results-table} may possibly be improved by considering additional information that is
not contained in our linear program.  

For instance, in 
Section \ref{tri-cross-sq} we were able to improve the bounds for
the triangle cross square by noting that the larger simplices in that
polytope must overlap.  Such considerations may
close the gap between upper and lower bounds for minimal
triangulations of simplotopes in several specific dimensions.  For
example, for $\Delta^2 \times \Delta^2 \times \Delta^1$, our lower
bound is 20 and the standard triangulation gives a construction of
size 30.  For $\Delta^2 \times \Delta^2 \times \Delta^2$, our lower
bound is 50, and by comparison the standard construction has size 90.

Also, the bounds that we used in our linear program for $V(s,t)$ only
relied on known results for the volumes of simplices in cubes (i.e.,
the Hadamard determinant problem);
however, bounds for $V(s,t)$ might perhaps be improved by restricting
attention to volumes of simplices in simplotopes (embedded in cubes).

Another direction that may improve our bounds
slightly is to consider corner simplices, as was done in
\cite{blsu05}.  The idea is that corner simplices in a cube have a large
number of exterior faces, and no other simplices have nearly as many.
There ought to be an analogous result for simplotopes.

Finally, we remark that our methods for products of segments and triangles may be generalized further to
consider products of other kinds of simplices, 
by more extensive bookkeeping of the interactions between simplices of various dimensions.

\subsection*{Acknowledgements}
The authors wish to thank Deborah Seacrest %and Francisco Santos 
for helpful comments.

\nocite{ramb02}
\nocite{lata82}
\nocite{blsu05}
\nocite{seacrest17}
\bibliography{bib}

\end{document}